\newif\iftechreport
\techreporttrue
\PassOptionsToPackage{table}{xcolor}

\documentclass{siamart}

\usepackage{dpr,fec}
\usepackage[noend]{algpseudocode}
\usepackage{xcolor}
\hypersetup{colorlinks}

\usepackage{rotating}
\usepackage{multirow}

\newcommand{\computesteps} {\textsc{Compute\_Steps}}
\newcommand{\computesigmas}{\textsc{Compute\_Sigma}}
\newcommand{\Fiteration}[1]{\textsc{F-iteration}#1}
\newcommand{\Viteration}[1]{\textsc{V-iteration}#1}

\newcommand{\deltaf}       {\delta^{f}}
\newcommand{\deltav}       {\delta^{v}}
\newcommand{\deltas}       {\delta^{s}}

\newcommand{\deltamax}     {\Delta}
\newcommand{\deltamaxv}    {\Delta^v}
\newcommand{\gammac}       {\gamma_{\rm c}}
\newcommand{\gammae}       {\gamma_{\rm e}}
\newcommand{\gammalambda}  {\gamma_{\lambda}}
\newcommand{\kappadelta}   {\kappa_{\delta}}
\newcommand{\kappafm}      {\kappa_{fm}}
\newcommand{\kappahs}      {\kappa_{hs}}
\newcommand{\kappaht}      {\kappa_{ht}}
\newcommand{\kappan}       {\kappa_{n}}
\newcommand{\kappantn}     {\kappa_{ntn}}
\newcommand{\kappantt}     {\kappa_{ntt}}
\newcommand{\kappap}       {\kappa_{p}}
\newcommand{\kapparho}     {\kappa_{\rho}}
\newcommand{\kappast}      {\kappa_{st}}
\newcommand{\kappavone}    {\kappa_{v1}}
\newcommand{\kappavm}      {\kappa_{vm}}
\newcommand{\kappavtwo}    {\kappa_{v2}}
\newcommand{\lambdaf}      {\lambda^{f}}
\newcommand{\lambdav}      {\lambda^{v}}
\newcommand{\sigmav}       {\sigma^{v}}
\newcommand{\sigmamin}     {\underline\sigma}
\newcommand{\sigmamax}     {\overline\sigma}
\newcommand{\thetafc}      {\theta_{fc}}
\newcommand{\Xiv}          {\Xi^v}
\newcommand{\xiv}          {\xi^v}

\newcommand{\gv}           {g^v}

\newcommand{\gp}           {g^p}
\newcommand{\Hv}           {H^v}
\newcommand{\mv}           {m^v}
\newcommand{\mf}           {m^f}
\newcommand{\ARC}          {\textsc{arc}}
\newcommand{\TRACE}        {\textsc{trace}}
\newcommand{\TrustFunnel}  {\textsc{TF}}
\newcommand{\TrustVonly}   {\textsc{TF-V-only}}
\newcommand{\contract}     {\textsc{contract}}
\newcommand{\yf}           {y^{f}}
\newcommand{\kSone}        {k_{\Scal_1}}
\newcommand{\kStwo}        {k_{\Scal_2}}
\newcommand{\kSthree}      {k_{\Scal_3}}
\newcommand{\NR}           {n^{R}}
\newcommand{\NN}           {n^{N}}
\newcommand{\tNull}        {t^{N}}

\newcommand{\AlgBreak}{\vspace*{-.7\baselineskip}\Statex\hspace*{\dimexpr-\algorithmicindent-2pt\relax}\rule{\textwidth}{0.4pt}}
\algrenewcommand{\algorithmiccomment}[1]{\hfill[#1]}

\title{Complexity Analysis of a Trust Funnel Algorithm for Equality Constrained Optimization\thanks{This material is based upon work supported by the U.S. Department of Energy, Office of Science, Office of Advanced Scientific Computing Research, Applied Mathematics, Early Career Research Program under contract number DE--SC0010615, as well as by the U.S. National Science Foundation under Grant No.~DMS--1319356.}}
\author{Frank~E.~Curtis\thanks{Department of Industrial and Systems Engineering, Lehigh University, Bethlehem, PA, USA.  E-mail: \email{frank.e.curtis@gmail.com}; \email{mos213@lehigh.edu}}
   \and Daniel~P.~Robinson\thanks{Department of Applied Mathematics and Statistics, Johns Hopkins University, Baltimore, MD, USA.  E-mail: \email{daniel.p.robinson@jhu.edu}}
   \and Mohammadreza Samadi\footnotemark[2]}
\date{\today}

\begin{document}

\maketitle

\begin{abstract}
  A method is proposed for solving equality constrained nonlinear optimization problems involving twice continuously differentiable functions.  The method employs a trust funnel approach consisting of two phases: a first phase to locate an $\epsilon$-feasible point and a second phase to seek optimality while maintaining at least $\epsilon$-feasibility.  A two-phase approach of this kind based on a cubic regularization methodology was recently proposed along with a supporting worst-case iteration complexity analysis.  Unfortunately, however, in that approach, the objective function is completely ignored in the first phase when $\epsilon$-feasibility is sought.  The main contribution of the method proposed in this paper is that the same worst-case iteration complexity is achieved, but with a first phase that also accounts for improvements in the objective function.  As such, the method typically requires fewer iterations in the second phase, as the results of numerical experiments demonstrate.
\end{abstract}



\section{Introduction}\label{sec.introduction}

The purpose of this paper is to propose a new method for solving equality constrained nonlinear optimization problems.  As is well known, such problems are important throughout science and engineering, arising in areas such as network flow optimization~\cite{hsu2002network,nygard2001dynamic}, optimal allocation with resource constraints~\cite{cui2004optimal,marti2004optimal}, maximum likelihood estimations with constraints~\cite{hathaway1985constrained}, and optimization with constraints defined by partial differential equations~\cite{biegler2007real,biros2005parallel,rees2010optimal}.

Contemporary methods for solving equality constrained optimization problems are predominantly based on ideas of sequential quadratic optimization (commonly known as SQP)~\cite{byrd2008inexact,CurJRW14,curtis2009matrix,GouR10a,GouR10b,GouR10c,morales2011sequential,NoceWrig06}.  The design of such methods remains an active area of research as algorithm developers aim to propose new methods that attain global convergence guarantees under weak assumptions about the problem functions.  Recently, however, researchers are being drawn to the idea of designing algorithms that also offer improved worst-case iteration complexity bounds.  This is due to the fact that, at least for convex optimization, algorithms designed with complexity bounds in mind have led to methods with improved practical performance.

For solving equality constrained optimization problems, a cubic regularization method is proposed in \cite{CartGoulToin13a} with an eye toward achieving good complexity properties.  This is a two-phase approach with a first phase that seeks an $\epsilon$-feasible point and a second phase that seeks optimality while maintaining $\epsilon$-feasibility.  The number of iterations that the method requires in the first phase to produce an $\epsilon$-feasible point is $\Ocal(\epsilon^{-3/2})$, a bound that is known to be optimal for unconstrained optimization \cite{CartGoulToin11c}.  The authors of \cite{CartGoulToin13a} then also propose a method for the second phase and analyze its complexity properties.  (For related work on cubic regularization methods for solving constrained optimization problems, see~\cite{CartGoulToin12,CartGoulToin14a}.)


Unfortunately, however, the method in \cite{CartGoulToin13a} represents a departure from the current state-of-the-art SQP methods that offer the best practical performance.  One of the main reasons for this is that contemporary SQP methods seek feasibility and optimality simultaneously.  By contrast, one of the main reasons that the approach from~\cite{CartGoulToin13a} does not offer practical benefits is that the first phase of the algorithm entirely ignores the objective function, meaning that numerous iterations might need to be performed before the objective function influences the trajectory of the algorithm.

The algorithm proposed in this paper can be considered a next step in the design of \emph{practical} algorithms for equality constrained optimization with good worst-case iteration complexity properties.  Ours is also a two-phase approach, but is closer to the SQP-type methods representing the state-of-the-art for solving equality constrained problems.  In particular, the first phase of our proposed approach follows a trust funnel methodology that locates an $\epsilon$-feasible point in $\Ocal(\epsilon^{-3/2})$ iterations \emph{while also attempting to yield improvements in the objective function}.  Borrowing ideas from the trust region method known as \TRACE{}~\cite{CurtRobiSama16}, we prove that our method attains the same worst-case iteration complexity bounds as those offered by~\cite{CartGoulToin13a}, and show with numerical experiments that consideration of the objective function in the first phase typically results in the second phase requiring fewer iterations.

\subsection{Organization}

In the remainder of this section,  we introduce notation that is used throughout the remainder of the paper and cover preliminary material on equality constrained nonlinear optimization.  In \S\ref{sec.phase1}, we motivate and describe our proposed ``phase 1'' method for locating an $\epsilon$-feasible point while also attempting to reduce the objective function.  An analysis of the convergence and worst-case iteration complexity of this phase 1 method is presented in \S\ref{sec.convergence}.  Strategies and corresponding convergence/complexity guarantees for ``phase 2'' are the subject of \S\ref{sec.phase2}, the results of numerical experiments are provided in \S\ref{sec.numerical}, and concluding remarks are given in \S\ref{sec.conclusion}.

\subsection{Notation}

Let $\R{}$ denote the set of real numbers (i.e., scalars), let $\R{}_{+}$ denote the set of nonnegative real numbers, let $\R{}_{++}$ denote the set of positive real numbers, and let $\N{} := \{1,2,\dots\}$ denote the set of natural numbers.  For any of these quantities, let a superscript $N \in \N{}$ be used to indicate the $N$-dimensional extension of the set---e.g., let $\R{N}$ denote the set of $N$-dimensional real vectors---and let a superscript $M \times N$ with $(M,N) \in \N{} \times \N{}$ be used to indicate the $M$-by-$N$-dimensional extension of the set---e.g., let $\R{M \times N}$ denote the set of $M$-by-$N$ real matrices.

A vector with all elements equal to 1 is denoted as $e$ and an identity matrix is denoted as $I$, where, in each case, the size of the quantity is determined by the context in which it appears.  With real symmetric matrices $A$ and $B$, let $A$ $\succ$ ($\succeq$) $B$ indicate that $A - B$ is positive definite (semidefinite); e.g., $A$ $\succ$ ($\succeq$) $0$ indicates that $A$ is positive definite (semidefinite).  
Given vectors $\{u,v\}\subset\Re^N$, let $u \perp v$ mean that $u_i v_i = 0$ for all $i\in\{1,2,\dots, N\}$.  Let $\|x\|$ denote the $2$-norm of a vector $x$.

\subsection{Preliminaries}

Given an objective function $f : \R{N} \to \R{}$ and constraint function $c : \R{N} \to \R{M}$, we study the equality constrained optimization problem
\bequation\label{prob.opt}
  \min_{x \in \R{N}}\ f(x)\ \ \st\ \ c(x) = 0.
\eequation
At the outset, let us state the following assumption about the problem functions.

\bassumption\label{ass.differentiability}
  The functions $f$ and $c$ are twice continuously differentiable.
\eassumption

In light of Assumption~\ref{ass.differentiability}, we define $g : \R{N} \to \R{N}$ as the gradient function of $f$, i.e., $g := \nabla f$, and define $J : \R{N} \to \R{M \times N}$ as the Jacobian function of $c$, i.e., $J := \nabla c^T$.  The function $c_i : \R{N} \to \R{}$ denotes the $i$th element of the function $c$.

Our proposed algorithm follows a local search strategy that merely aims to compute a first-order stationary point for problem~\eqref{prob.opt}.  Defining the Lagrangian $\Lcal : \R{N} \times \R{M} \to \R{}$ as given by $\Lcal(x,y) = f(x) + y^Tc(x)$, a first-order stationary point $(x,y)$ is one that satisfies $0 = \nabla_x L(x,y) \equiv g(x) + J(x)^Ty$ and $0 = \nabla_y L(x,y) \equiv c(x)$.

Our proposed technique for solving problem~\eqref{prob.opt} is iterative, generating, amongst other quantities, a sequence of iterates $\{x_k\}$ indexed by $k \in \N{}$.  For ease of exposition, we also apply an iteration index subscript for function and other quantities corresponding to the $k$th iteration; e.g., we write $f_k$ to denote $f(x_k)$.

\section{Phase 1: Obtaining Approximate Feasibility}\label{sec.phase1}

The goal of phase 1 is to obtain an iterate that is (approximately) feasible.  This can, of course, be accomplished by employing an algorithm that focuses exclusively on minimizing a measure of constraint violation.  However, we find this idea to be unsatisfactory since such an approach would entirely ignore the objective function.  Alternatively, in this section, we present a trust funnel algorithm with good complexity properties for obtaining (approximate) feasibility that attempts to simultaneously reduce the objective function, as is commonly done in contemporary nonlinear optimization algorithms.

\subsection{Step computation}

Similar to other trust funnel algorithms \cite{GoulToin10,CurtGoulRobiToin16}, our algorithm employs a step-decomposition approach wherein each trial step is composed of a \emph{normal step} aimed at reducing constraint violation (i.e., infeasibility) and a \emph{tangential step} aimed at reducing the objective function.  The algorithm then uses computed information, such as the reductions that the trial step yields in models of the constraint violation and objective function, to determine which of two types of criteria should be used for accepting or rejecting the trial step.  To ensure that sufficient priority is given to obtaining (approximate) feasibility, an upper bound on a constraint violation measure is initialized, maintained, and subsequently driven toward zero as improvements toward feasibility are obtained.  The algorithm might also nullify the tangential component of a trial step, even after it is computed, if it is deemed too harmful in the algorithm's pursuit toward (approximate) feasibility.  In this subsection, the details of our approach for computing a trial step are described.

\subsubsection{Normal step}

The purpose of the normal step is to reduce infeasibility.  The measure of infeasibility that we employ is $v : \R{N} \to \R{}$ defined by
\bequation \label{def.v}
  v(x) = \thalf \|c(x)\|^2.
\eequation
At an iterate $x_k$, the normal step $n_k$ is defined as a minimizer of a second-order Taylor series approximation of $v$ at $x_k$ subject to a trust region constraint, i.e.,
\bequation\label{prob.normal_step}
  n_k \in \arg \min_{n \in \R{N}} \mv_k(n) \ \ \st \ \ \|n\| \leq \deltav_k,
\eequation
where the scalar $\deltav_k \in (0,\infty)$ is the trust region radius and the model of the constraint violation measure at $x_k$ is $\mv_k : \R{N} \to \R{}$ defined by
\begin{align}
  \mv_k(n) =&\ v_k + {\gv_k}^T  n + \thalf n^T \Hv_k n
  \ \ \ \text{with} \ \ \ 
  \gv_k   := \nabla v(x_k) = J_k^T c_k, \label{def.gv} \\ \text{and}\ \ 
  \Hv_k   :=&\ \nabla^2 v(x_k) = J_k^T J_k + \sum_{i=1}^M c_i(x_k) \nabla^2 c_i(x_k). \label{def.Hv}
\end{align}
For any $(x_k, \deltav_k) \in \R{N} \times \R{}_+$, a globally optimal solution to \eqref{prob.normal_step} exists~\cite[Corollary~7.2.2]{ConGT00a} and $n_k$ has a corresponding dual variable $\lambdav_k \in \R{}_+$ such that
\bsubequations\label{eq.kkt_normal}
  \begin{align}
    \gv_k + (\Hv_k + \lambdav_k I) n_k &= 0, \label{eq.kkt-normal.dual} \\
    \Hv_k + \lambdav_k I  &\succeq 0, \label{eq.kkt-normal.psd} \\ \text{and}\ \ 
    0 \leq \lambdav_k \perp (\deltav_k - \|n_k\|) &\geq 0. \label{eq.kkt-normal.comp}
  \end{align}
\esubequations

In a standard trust region strategy, a trust region radius is given at the beginning of an iteration, which \emph{explicitly} determines the primal-dual solution of the subproblem.  Our method, on the other hand, might instead make use of a normal step that is derived as a solution of~\eqref{prob.normal_step} where the trust region radius is defined \emph{implicitly} by a given dual variable $\lambdav_k$.  In particular, given $\lambdav_k \in [0,\infty)$ that is strictly larger than the negative of the leftmost eigenvalue of $\Hv_k$, our algorithm might compute $n_k$ from
\bequation\label{prob.Q_v}
  \Qcal_k^v(\lambdav_k) :\ \ \min_{n\in\R{N}}\ v_k + {\gv_k}^T n + \thalf n^T(\Hv_k + \lambdav_k I)n.
\eequation
The unique solution to~\eqref{prob.Q_v}, call it $n_k(\lambdav_k)$, is the solution of the nonsingular linear system $(\Hv_k + \lambdav_k I)n = -\gv_k$, and is the global solution of~\eqref{prob.normal_step} for $\deltav_k = \|n_k(\lambdav_k)\|$.

\subsubsection{Tangential step}

The purpose of the tangential step is to reduce the objective function.  Specifically, when requested by the algorithm, the tangential step~$t_k$ is defined as a minimizer of a quadratic model of the objective function in the null space of the constraint Jacobian subject to a trust region constraint, i.e.,
\bequation\label{prob.tangential_step}
  t_k \in \arg \min_{t \in \R{N}} \mf_k (n_k+t)\ \ \st \ \ J_kt = 0\ \ \text{and}\ \ \|n_k+t\| \leq \deltas_k,
\eequation
where $\deltas_k \in (0,\infty)$ is a trust region radius and, with some symmetric $H_k \in \R{N\times N}$, the objective function model $\mf_k : \R{N} \to \R{}$ is defined by
\bequation\label{def.mf}
  \mf_k(s) = f_k + g_k^T s + \thalf s^T H_k s.
\eequation
Following other trust funnel strategies, one desires $\deltas_k$ to be set such that the trust region describes the area in which the models of the constraint violation and objective function are accurate.  In particular, with a trust region radius $\deltaf_k \in (0,\infty)$ for the objective function, our algorithm employs, for some scalar $\kappadelta \in (1,\infty)$, the value
\bequation\label{def.delta_s}
  \deltas_k := \min \{\kappadelta \deltav_k, \deltaf_k\}.
\eequation
Due to this choice of trust region radius, it is deemed not worthwhile to compute a nonzero tangential step if the feasible region of \eqref{prob.tangential_step} is small.  Specifically, our algorithm only computes a nonzero $t_k$ when $\|n_k\| \leq \kappan \deltas_k$ for some $\kappan \in (0,1)$. In addition, it only makes sense to compute a tangential step when reasonable progress in reducing~$f$ in the null space of $J_k$ can be expected.  To predict the potential progress, we define 
\bequation\label{def.gp}
  \gp_k := Z_kZ_k^T(g_k+H_k n_k),
\eequation
where the columns of $Z_k$ form an orthonormal basis for $\Null(J_k)$.  If $\|\gp_k\| < \kappap \|g^v_k\|$ for some $\kappap \in(0,\infty)$, then computing a tangential step is not worthwhile and we simply set the primal-dual solution (estimate) for \eqref{prob.tangential_step} to zero.

For any $(x_k,\deltas_k,H_k) \in \R{N} \times \R{}_+ \times \R{N \times N}$, a globally optimal solution to \eqref{prob.tangential_step} exists~\cite[Corollary~7.2.2]{ConGT00a} and $t_k$ has corresponding dual variables $\yf_k \in \R{M}$ and $\lambdaf_k \in \R{}_+$ (for the null space and trust region constraints, respectively) such that
\bsubequations\label{eq.kkt_tangent}
  \begin{align}
    \bbmatrix H_k + \lambdaf_k I & J_k^T \\ J_k & 0 \ebmatrix \bbmatrix t_k \\ \yf_k \ebmatrix &= - \bbmatrix g_k + (H_k + \lambdaf_k I) n_k \\ 0 \ebmatrix, \label{eq.kkt_tangent.dual} \\
    Z_k^TH_kZ_k + \lambdaf_k I &\succeq 0, \label{eq.kkt_tangent.psd} \\ \text{and}\ \ 
    0 \leq \lambdaf_k \perp ({\deltas_k} - \|n_k + t_k\|) &\geq 0. \label{eq.kkt_tangent.comp}
  \end{align}
\esubequations
Similarly as for the normal step computation, our algorithm potentially computes $t_k$ not as a solution of \eqref{prob.tangential_step} for a given $\deltas_k$, but as a solution of a regularized subproblem for a given dual variable for the trust region constraint.  Specifically, for $\lambdaf_k \in [0,\infty)$ that is strictly larger than the negative of the leftmost eigenvalue of $Z_k^TH_kZ_k$, our algorithm might solve the following subproblem for the tangential step:
\bequation\label{prob.Q_f}
  \Qcal_k^f(\lambdaf_k) :\ \ \min_{t \in \R{N}}\ (g_k + (H_k+\lambdaf_k I) n_k)^T t + \thalf t^T(H_k+\lambdaf_k I) t\ \ \st \ \ J_k t = 0.
\eequation
The unique solution $t_k(\lambdaf_k)$ of~\eqref{prob.Q_f} is a global solution of~\eqref{prob.tangential_step} for $\deltas_k = \|n_k + t_k(\lambdaf_k)\|$.

There are situations in which our algorithm discards a computed tangential step after one is computed, i.e., situations when the algorithm resets $t_k \gets 0$.  Specifically, this occurs when any of the following conditions fails to hold:
\bsubequations\label{eq.tangential_step_conditions}
  \begin{align}
    \mv_k(0) - \mv_k(n_k + t_k) &\geq \kappavm \(\mv_k(0) - \mv_k(n_k)\) &&\text{for some}\ \ \kappavm \in (0,1); \label{eq.mv_red_t} \\
    \|n_k + t_k\| &\geq \kappantn \|n_k\| &&\text{for some}\ \ \kappantn \in (0,1); \label{eq.s_vs_n} \\
    \|\Hv_kt_k\| &\leq \kappaht \|n_k+t_k\|^2 &&\text{for some}\ \ \kappaht \in (0,\infty). \label{eq.Ht_vs_s}
  \end{align}
\esubequations
The first of these conditions requires that the reduction in the constraint violation model for the full step $s_k := n_k + t_k$ is sufficiently large with respect to that obtained by the normal step; the second requires that the full step is sufficiently large in norm compared to the normal step; and the third requires that the action of the tangential step on the Hessian of the constraint violation model is not too large compared to the squared norm of the full step.  It is worthwhile to mention that all of these conditions are satisfied automatically when $\Hv_k = J_k^TJ_k$ (recall \eqref{def.Hv}), which occurs, e.g., when $c$ is affine.  However, since this does not hold in general, our algorithm requires these conditions explicitly, or else resets the tangential step to zero (which satisfies \eqref{eq.tangential_step_conditions}).

\subsection{Step acceptance}

After computing a normal step $n_k$ and potentially a tangential step $t_k$, the algorithm determines whether to accept the full step $s_k := n_k + t_k$.  The strategy that it employs is based on first using the obtained reductions in the models of constraint violation and the objective, as well as other related quantities, to determine what should be the main goal of the iteration: reducing constraint violation or the objective function.  Since the primary goal of phase~1 is to obtain (approximate) feasibility, the algorithm has a preference toward reducing constraint violation unless the potential reduction in the objective function is particularly compelling.  Specifically, if the following conditions hold, indicating good potential progress in reducing the objective, then the algorithm performs an \Fiteration{} (see \S\ref{sec.f-iteration}):
\bsubequations\label{eq.f-iteration-conditions}
  \begin{align}
    t_k \neq 0\ \ \text{with}\ \ \|t_k\| &\geq \kappast\|s_k\| && \text{for some}\ \kappast \in (0,1), \label{eq.f-iter_t_vs_s} \\
    \mf_k(0) - \mf_k(s_k) &\geq \kappafm (\mf_k(n_k) - \mf_k(s_k)), && \text{for some}\ \kappafm \in (0,1), \label{eq.f-iter_mf_red_t} \\
    v(x_k + s_k) &\leq \vmax_k - \kapparho \|s_k\|^3 && \text{for some}\ \kapparho \in (0,1),  \label{eq.f-iter_vmax_rule} \\
    n_k^T t_k &\geq -\thalf \kappantt \|t_k\|^2 && \text{for some}\ \kappantt \in (0,1), \label{eq.f-iter_nt_vs_t} \\
    \lambdav_k &\leq \sigmav_k\|n_k\| && \text{and} \label{eq.f-iter_sigma_v} \\
    \left\|(H_k-\nabla^2 f(x_k))s_k\right\| &\leq \kappahs \|s_k\|^2 && \text{for some}\ \kappahs \in (0,\infty). \label{eq.f-iter_Hs_vs_s}
  \end{align}
\esubequations
Conditions \eqref{eq.f-iter_t_vs_s}--\eqref{eq.f-iter_vmax_rule} are similar to those employed in other trust funnel algorithms, except that \eqref{eq.f-iter_t_vs_s} and \eqref{eq.f-iter_vmax_rule} are stronger (than the common, weaker requirements that $t_k \neq 0$ and $v(x_k+s_k) \leq \vmax_k$).  Employed here is a scalar sequence $\{\vmax_k\}$ updated dynamically by the algorithm that represents an upper bound on constraint violation; for this sequence, the algorithm ensures (see Lemma~\ref{lem.vmax}) that $v_k \leq \vmax_k$ and $\vmax_{k+1} \leq \vmax_k$ for all $k \in \N{}$.  Condition \eqref{eq.f-iter_nt_vs_t} ensures that, for an \Fiteration{}, the inner product between the normal and tangential steps is not too negative (or else the tangential step might undo too much of the progress toward feasibility offered by the normal step).  Finally, conditions \eqref{eq.f-iter_sigma_v} and \eqref{eq.f-iter_Hs_vs_s} are essential for achieving good complexity properties, requiring that any \Fiteration{} involves a normal step that is sufficiently large compared to the Lagrange multiplier for the trust region constraint and that the action of the full step on $H_k$ does not differ too much from its action on $\nabla^2 f(x_k)$.  If any condition in~\eqref{eq.f-iteration-conditions} does not hold, then a \Viteration{} is perfomed (see \S\ref{sec.v-iteration}).

Viewing \eqref{eq.f-iter_Hs_vs_s}, it is worthwhile to reflect on the choice of $H_k$ in the algorithm.  With the attainment of optimality (not only feasibility) in mind, standard practice would suggest that it is desirable to choose $H_k$ as the Hessian of the Lagrangian of problem~\eqref{prob.opt} for some multiplier vector $y_k \in \R{M}$.  This multiplier vector could be obtained, e.g., as the \emph{QP multipliers} from some previous iteration or \emph{least squares multipliers} using current derivative information.  For our purposes of obtaining good complexity properties for phase 1, we do not require a particular choice of $H_k$, but this discussion and \eqref{eq.f-iter_Hs_vs_s} do offer some guidance.  Specifically, one might choose $H_k$ as an approximation of the Hessian of the Lagrangian, potentially with the magnitude of the multiplier vector restricted in such a way that, after the full step is computed, the action of it on $H_k$ will not differ too much with its action on $\nabla^2 f(x_k)$.

\subsubsection{\Fiteration{}}\label{sec.f-iteration}

If~\eqref{eq.f-iteration-conditions} holds, then we determine that the $k$th iteration is an \Fiteration.  In this case, we begin by calculating the quantity
\bequation\label{def.rho_f}
  \rho_k^f \gets (f_k - f(x_k + s_k))/\|s_k\|^3,
\eequation
which is a measure of the decrease in $f$.  Using this quantity, acceptance or rejection of the step and the rules for updating the trust region radius are similar as in~\cite{CurtRobiSama16}.  As for updating the trust funnel radius, rather than the update in~\cite[Algorithm 2.1]{GoulToin10}, we require a modified update to obtain our complexity result; in particular, we use
\bequation\label{eq.vmax_update_f}
  \vmax_{k+1} \gets \min\{\max\{\kappavone \vmax_k,\vmax_k - \kapparho \|s_k\|^3\}, v_{k+1} + \kappavtwo(\vmax_k - v_{k+1})\}
\eequation
for some $\{\kappavone,\kappavtwo\}\subset (0,1)$.

\subsubsection{\Viteration{}}\label{sec.v-iteration}

When any one of the conditions in~\eqref{eq.f-iteration-conditions} does not hold, the $k$th iteration is a \Viteration{}, during which the main focus is to decrease the measure of infeasibility $v$.  In this case, we calculate
\bequation\label{def.rho_v}
  \rho_k^v \gets (v(x_k) - v(x_k + s_k))/\|s_k\|^3,
\eequation
which provides a measure of the decrease in constraint violation.  The rules for accepting or rejecting the trial step and for updating the trust region radius are the same as those in~\cite{CurtRobiSama16}.  One addition is that during a successful \Viteration{}, the trust funnel radius is updated, using the same constants as in \eqref{eq.vmax_update_f}, as
\bequation \label{eq.vmax_update_v}
  \vmax_{k+1} \gets \min\{\max\{\kappavone \vmax_k,v_{k+1} + \kappavtwo(v_k - v_{k+1})\}, v_{k+1} + \kappavtwo(\vmax_k - v_{k+1})\}.
\eequation

\subsection{Algorithm statement}

Our complete algorithm for finding an (approximately) feasible point can now be stated as Algorithm~\ref{alg.phase1} on page~\pageref{alg.phase1}, which in turn calls the \Fiteration{} subroutine stated as Algorithm~\ref{alg.Fiter} on page~\pageref{alg.Fiter} and the \Viteration{} subroutine stated as Algorithm~\ref{alg.viter} on page~\pageref{alg.viter}.

\balgorithm[ht]
  \small
  \caption{Trust Funnel Algorithm for Phase 1}
  \label{alg.phase1}
  \balgorithmic[1]
    \smallskip
    \Require $\{\kappan,\kappavm,\kappantn,\kapparho,\kappafm,\kappast,\kappantt,\kappavone,\kappavtwo,\gammac\}\subset(0,1)$,
    \Statex \hspace{26pt} $\{\kappap,\kappaht,\kappahs,\epsilon,\underline{\sigma}\} \subset (0,\infty)$, $\{\kappadelta,\gammae,\gammalambda\} \in (1,\infty)$, and $\overline{\sigma} \in [\underline{\sigma},\infty)$;
    \Statex \hspace{26pt} \Fiteration\ (Algorithm~\ref{alg.Fiter}, page~\pageref{alg.Fiter}) and \Viteration\ (Algorithm~\ref{alg.viter}, page~\pageref{alg.viter})
    \smallskip
    \AlgBreak
    \Procedure{Trust\_Funnel}{}
      \State choose $x_0 \in \R{N}$, $\vmax_0 \in [\max\{1,v_0\},\infty)$, and $\sigmav_0 \in [\sigmamin,\sigmamax]$
      \State choose $\{\deltav_0,\deltamaxv_0,\deltaf_0\} \subset (0,\infty)$ such that $\deltav_0 \leq \deltamaxv_0$ 
      \For {$k \in \N{}$}
        \If {$\|g^v_k\| \leq \epsilon$}
          \State \label{step.terminate} \Return $x_k$ 
        \EndIf
        \State $(n_k,t_k,\lambdav_k,\lambdaf_k) \gets \computesteps(x_k,\deltav_k,\deltas_k)$
        \State set $s_k \gets n_k + t_k$ \label{step.phase1.s_def}
        \State set $\sigmav_k \gets \computesigmas(n_k,\lambdav_k,\sigmav_{k-1},\rho_{k-1}^v)$
        \If{\eqref{eq.f-iteration-conditions} is satisfied} \label{step.phase1.check_type}
          \State set $\rho_k^f$ by \eqref{def.rho_f}
          \State $(x_{k+1},\vmax_{k+1},\deltaf_{k+1}) \gets \Fiteration(x_k,n_k,s_k,\vmax_k,\deltaf_k,\lambdaf_k,\rho_k^f)$
          \State set $\deltav_{k+1} \gets \deltav_k$, $\deltamaxv_{k+1} \gets \deltamaxv_k$, and $\rho_k^v \gets \infty$ \label{step.phase1.f-iter-updates}
        \Else
          \State set $\rho_k^v$ by \eqref{def.rho_v} \label{step.phase1.v-iter-rho}
          \State $(x_{k+1},\vmax_{k+1},\deltav_{k+1},\deltamaxv_{k+1}) \gets \Viteration(x_k,n_k,s_k,\vmax_k,\deltav_k,\deltamaxv_k,\lambdav_k,\sigmav_k,\rho_k^v)$
          \State set $\deltaf_{k+1} \gets \deltaf_k$ and $\rho_k^f \gets \infty$ \label{step.phase1.v-iter-updates}
        \EndIf
      \EndFor
    \EndProcedure
    \smallskip
    \AlgBreak
    \Procedure{\computesteps}{$x_k,\deltav_k,\deltas_k$}
      \State set $(n_k,\lambdav_k)$ as a primal-dual solution to \eqref{prob.normal_step}
      \State set $(t_k,\lambdaf_k) \gets (0,0)$
      \If {$\|n_k\| \leq \kappan \deltas_k$ and $\|\gp_k\| \geq \kappap \|g^v_k\|$}\label{step.phase1.tangential_conditions}
        \State set $(t_k,\yf_k,\lambdaf_k)$ as a primal-dual solution to~\eqref{prob.tangential_step}
        \State \textbf{if} any condition in~\eqref{eq.tangential_step_conditions} fails to hold \textbf{then} set $(t_k,\lambdaf_k) \gets (0,0)$ \label{step.phase1.tangential_reset}
      \EndIf
      \State \Return $(n_k,t_k,\lambdav_k,\lambdaf_k)$
    \EndProcedure
    \smallskip
    \AlgBreak
    \Procedure{\computesigmas}{$n_k,\lambdav_k,\sigmav_{k-1},\rho_{k-1}^v$}
    \If{iteration $(k-1)$ was an \Fiteration}
      \State set $\sigmav_k \gets \sigmav_{k-1}$ \label{step.phase1.f-iter_sigmav}
    \Else
      \State \textbf{if} $\rho_{k-1}^v < \kapparho$ \textbf{then} set $\sigmav_{k} \gets \max\{\sigmav_{k-1},\lambdav_{k}/\|n_{k}\|\}$ \textbf{else} set $\sigmav_{k} \gets \sigmav_{k-1}$ \label{step.phase1.v-iter_sigmav}
    \EndIf
    \State \Return $\sigmav_k$
  \EndProcedure
  \ealgorithmic
\ealgorithm

\balgorithm[ht]
  \small
  \caption{\Fiteration{} subroutine}
  \label{alg.Fiter}
  \balgorithmic[1]
    \smallskip
    \Procedure{\Fiteration}{$x_k,n_k,s_k,\vmax_k,\deltaf_k,\lambdaf_k,\rho_k^f$}
      \If{$\rho_k^f \geq \kapparho$} \label{step.accept-Fiter} \Comment{accept step}
        \State set $x_{k+1} \gets x_k+s_k$ \label{step.accept_update-Fiter}
        \State set $\vmax_{k+1}$ according to~\eqref{eq.vmax_update_f} \label{step.vmaxupdate-Fiter}
        \State set $\deltaf_{k+1} \gets \max\{\deltaf_k,\gammae \|s_k\|\}$ \label{step.delta_success-Fiter}
      \Else\ (i.e., \textbf{if} $\rho_k^f < \kapparho$) \label{step.contract-Fiter} \Comment{contract trust region}
        \State set $x_{k+1} \gets x_k$
        \State set $\vmax_{k+1} \gets \vmax_k$
        \State set $\deltaf_{k+1} \gets \text{F-}\contract(n_k,s_k,\deltaf_k,\lambdaf_k)$ \label{step.delta_contract-Fiter}
      \EndIf
      \State \Return $(x_{k+1},\vmax_{k+1},\deltaf_{k+1})$
    \EndProcedure
    \smallskip
    \AlgBreak
    \Procedure{F-\contract}{$n_k,s_k,\deltaf_k,\lambdaf_k$}
      \If{$\lambdaf_k < \sigmamin \|s_k\|$} \label{step.lamcheck-Fiter}
        \State set $\lambdaf > \lambdaf_k$ so the solution $t(\lambdaf)$ of $\Qcal_k^f(\lambdaf)$ yields $\sigmamin \leq \lambdaf/\|n_k+t(\lambdaf)\|$ \label{step.linsysmagic-Fiter}
        \State \Return $\deltaf_{k+1} \gets \|n_k+t(\lambdaf)\|$ \label{step.deltamagic-Fiter}
      \Else\ (i.e., \textbf{if} $\lambdaf_k \geq \sigmamin \|s_k\|$) \label{step.else-Fiter}
        \State \Return $\deltaf_{k+1} \gets \gammac\|s_k\|$ \label{step.deltanormal-Fiter}
      \EndIf
    \EndProcedure
  \ealgorithmic
\ealgorithm

\balgorithm[ht]
  \small
  \caption{\Viteration{} subroutine}
  \label{alg.viter}
  \balgorithmic[1]
    \smallskip
    \Procedure{\Viteration{}}{$x_k,n_k,s_k,\vmax_k,\deltav_k,\deltamaxv_k,\lambdav_k,\sigmav_k,\rho_k^v$}
      \If{$\rho_k^v \geq \kapparho$ and either $\lambdav_k \leq \sigmav_k \|n_k\|$ or $\|n_k\| = \deltamaxv_k$} \label{step.accept-Viter} \Comment{accept step}
        \State set $x_{k+1} \gets x_k+s_k$ \label{step.accept_update-Viter}
        \State set $\vmax_{k+1}$ according to~\eqref{eq.vmax_update_v}
        \State set $\deltamaxv_{k+1} \gets \max\{\deltamaxv_k,\gammae\|n_k\|\}$ \label{step.deltamaxA-Viter}
        \State set $\deltav_{k+1} \gets \min\{\deltamaxv_{k+1},\max\{\deltav_k, \gammae\|n_k\|\}\}$ \label{step.delta_success-Viter}
      \ElsIf{$\rho_k^v < \kapparho$} \label{step.contract-Viter} \Comment{contract trust region}
        \State set $x_{k+1} \gets x_k$
        \State set $\vmax_{k+1} \gets \vmax_k$
        \State set $\deltamaxv_{k+1} \gets \deltamaxv_k$ \label{step.deltamaxC-Viter}
        \State set $\deltav_{k+1} \gets \text{V-}\contract(n_k,s_k,\deltav_k,\lambdav_k)$ \label{step.delta_contract-Viter}
      \Else\ (i.e., \textbf{if} $\rho_k^v \geq \kapparho$, $\lambdav_k > \sigmav_k \|n_k\|$, and $\|n_k\| < \deltamaxv_k$) \label{step.expand-Viter} \Comment{expand trust region}
        \State set $x_{k+1} \gets x_k$
        \State set $\vmax_{k+1} \gets \vmax_k$
        \State set $\deltamaxv_{k+1} \gets \deltamaxv_k$ \label{step.deltamaxE-Viter}
        \State set $\deltav_{k+1} \gets \min\{\deltamaxv_{k+1},\lambdav_k/\sigmav_k\}$ \label{step.delta_expand-Viter}
      \EndIf
      \State \Return $(x_{k+1},\vmax_{k+1},\deltav_{k+1},\deltamaxv_{k+1})$
    \EndProcedure
    \smallskip
    \AlgBreak
    \Procedure{V-\contract}{$n_k,s_k,\deltav_k,\lambdav_k$}
      \If{$\lambdav_k < \sigmamin \|n_k\|$} \label{step.lamcheck-Viter}
        \State set $\hat{\lambda}^v \gets \lambdav_k + (\sigmamin \|\gv_k\|)^{1/2}$ \label{step.lambdahat-Viter}
        \State set $\lambdav \gets \hat{\lambda}^v$ \label{step.lambda0-Viter} \label{step.lambdav1}
        \State set $n(\lambdav)$ as the solution of $\Qcal_k^v(\lambdav)$ \label{step.find-n}
        \If{$\lambdav / \|n(\lambdav)\| \leq \sigmamax$} \label{step.checkratio-Viter}
          \State \Return $\deltav_{k+1} \gets \|n(\lambdav)\|$ \label{step.regularized1-Viter}
        \Else
          \State set $\lambdav \in (\lambdav_k,\hat{\lambda}^v)$ so the solution $n(\lambdav)$ of $\Qcal_k^v(\lambdav)$ yields $\sigmamin \leq \lambdav/\|n(\lambdav)\| \leq \sigmamax$ \label{step.linsysmagic-Viter}
          \State \Return $\deltav_{k+1} \gets \|n(\lambdav)\|$ \label{step.deltamagic-Viter}
        \EndIf
      \Else\ (i.e., \textbf{if} $\lambdav_k \geq \sigmamin \|n_k\|$) \label{step.else-Viter}
        \State \label{step.lambda-Viter} set $\lambdav \gets \gammalambda\lambdav_k$
        \State \label{step.linsys-Viter} set $n(\lambdav)$ as the solution of $\Qcal^v_k(\lambdav)$
        \If{$\|n(\lambdav)\| \geq \gammac\|n_k\|$} \label{eq.sabotage}
          \State \Return $\deltav_{k+1} \gets \|n(\lambdav)\|$ \label{step.delta2-Viter}
        \Else
          \State \Return $\deltav_{k+1} \gets \gammac\|n_k\|$ \label{step.deltanormal-Viter}
        \EndIf
      \EndIf
    \EndProcedure
  \ealgorithmic
\ealgorithm

\section{Convergence and Complexity Analyses for Phase 1}\label{sec.convergence}

The analyses that we present require the following assumption related to the iterate sequence.

\bassumption\label{ass.compactness}
  The sequence of iterates $\{x_k\}$ is contained in a compact set.  In addition, the sequence $\{\|H_k\|\}$ is bounded over $k \in \N{}$.
\eassumption

Our analysis makes extensive use of the following mutually exclusive and exhaustive subsets of the iteration index sequence generated by Algorithm~\ref{alg.phase1}:
\bequalin
  \Ical &:= \{k \in \N{}  : \|g^v_k\| > \epsilon\}, \\
  \Fcal &:= \{k \in \Ical : \text{iteration $k$ is an \Fiteration{}}\}, \\ \text{and}\ \ 
  \Vcal &:= \{k \in \Ical : \text{iteration $k$ is a \Viteration{}}\}.
\eequalin
It will also be convenient to define the index set of iterations for which tangential steps are computed and not reset to zero by our method:
\bequalin
  \Ical^t 
    &:=           \{k \in \Ical : t_k \neq 0 \ \text{when Step~\ref{step.phase1.s_def} of Algorithm~\ref{alg.phase1} is reached} \} \\
    &\phantom{:}= \{k \in \Ical : \text{\ Step~\ref{step.phase1.tangential_reset} of Algorithm~\ref{alg.phase1} is reached and all conditions in \eqref{eq.tangential_step_conditions} hold}\}.
\eequalin

\subsection{Convergence analysis for phase 1}

The goal of our convergence analysis is to prove that Algorithm~\ref{alg.phase1} terminates finitely, i.e., $|\Ical| < \infty$.  Our analysis to prove this fact requires a refined examination of the subsets $\Fcal$ and $\Vcal$ of~$\Ical$.  For these purposes, we define disjoint subsets of $\Fcal$ as
\bequationn
  \Scal^f := \{k \in \Fcal : \rho_k^f \geq \kapparho\}\ \ \text{and}\ \ \Ccal^f := \{k \in \Fcal : \rho_k^f < \kapparho\},
\eequationn
and disjoint subsets of $\Vcal$ as
\begin{align*}  
  \Scal^v &:= \{k \in \Vcal : \rho_k^v \geq \kapparho \text{ and either } \lambdav_k \leq \sigmav_k \|n_k\| \text{ or }\|n_k\| = \deltamaxv_k\}, \\
  \Ccal^v &:= \{k \in \Vcal : \rho_k^v < \kapparho\}, \\ \text{and}\ \ 
  \Ecal^v &:= \{k \in \Vcal : k \notin \Scal^v \cup \Ccal^v\}.
\end{align*}
We further partition the set $\Scal^v$ into the two disjoint subsets
\bequationn
  \Scal^v_{\deltamax} := \{k \in \Scal^v : \|n_k\| = \deltamaxv_k\}\ \ \text{and}\ \ \Scal^v_{\sigma} := \{k \in \Scal^v : k \notin \Scal^v_{\deltamax}\}.
\eequationn
Finally, for convenience, we also define the unions
\bequationn
  \Scal := \{k \in \Ical : k \in \Scal^f \cup \Scal^v\}\ \ \text{and}\ \ \Ccal := \{k \in \Ical : k \in \Ccal^f \cup \Ccal^v\}.
\eequationn
Due to the updates for the primal iterate and/or trust region radii in the algorithm, we often refer to iterations with indices in $\Scal$ as successful steps, those with indices in~$\Ccal$ as contractions, and those with indices in $\Ecal^v$ as expansions.

Basic relationships between all of these sets are summarized in our first lemma.

\blemma
  The following relationships hold:
  \bitemize
    \item[(i)]   $\Fcal \cap \Vcal = \emptyset$ and $\Fcal \cup \Vcal = \Ical$;
    \item[(ii)]  $\Scal^f \cap \Ccal^f = \emptyset$ and $\Scal^f \cup \Ccal^f = \Fcal$;
    \item[(iii)] $\Scal^v$, $\Ccal^v$, and $\Ecal^v$ are mutually disjoint and $\Scal^v \cup \Ccal^v \cup \Ecal^v = \Vcal$; and
    \item[(iv)]  if $k \in \Ical \setminus \Ical^t$, then $k \in \Vcal$.
  \eitemize
\elemma
\bproof
  The fact that $\Fcal \cap \Vcal = \emptyset$ follows from the two cases resulting from the conditional statement in Step~\ref{step.phase1.check_type} of Algorithm~\ref{alg.phase1}.  The rest of part (i), part (ii), and part (iii) follow from the definitions of the relevant sets. Part (iv) can be seen to hold as follows.  If $k\in\Ical \setminus \Ical^t$, then $t_k = 0$ so that~\eqref{eq.f-iter_t_vs_s} does not hold. It now follows from the logic in Algorithm~\ref{alg.phase1} that $k\in\Vcal$ as claimed.
\eproof

The results in the next lemma are consequences of Assumptions~\ref{ass.differentiability} and \ref{ass.compactness}.

\blemma\label{lem.derivative_bounds}
  The following hold:
  \bitemize
    \item[(i)]   there exists $\thetafc \in (1,\infty)$ so $\max\{\|g_k\|,\|c_k\|,\|J_k\|,\|\Hv_k\|\} \leq \thetafc$ for all $k \in \Ical$;
    \item[(ii)]  $\|\gv_k\| \equiv \|J_k^T c_k\| \leq \thetafc \|c_k\|$ for all $k\in\Ical$; and
    \item[(iii)] $\gv : \R{N} \to \R{N}$ defined by $\gv(x) = J(x)^T c(x)$ (recall \eqref{def.gv}) is Lipschitz continuous with Lipschitz constant $\gv_{Lip} > 0$ over an open set containing $\{x_k\}$.
  \eitemize
\elemma
\bproof
  Part (i) follows from Assumptions~\ref{ass.differentiability} and~\ref{ass.compactness}.  Part~(ii) follows since, by the Cauchy--Schwarz inequality, $\|J_k^T c_k\| \leq \|J_k\| \|c_k\| \leq \thetafc \|c_k\|$.  Part (iii) follows since the first derivative of $g^v$ is uniformly bounded under Assumptions~\ref{ass.differentiability} and~\ref{ass.compactness}.
\eproof

We now summarize properties associated with the normal and tangential steps.

\blemma\label{lem.nonzero_steps}
  The following hold for all $k\in\Ical$:
  \bitemize
    \item[(i)]  $n_k \neq 0$ and $s_k \neq 0$; and
    \item[(ii)] in Step~\ref{step.phase1.s_def} of Algorithm~\ref{alg.phase1}, the vector $t_k$ satisfies~\eqref{eq.tangential_step_conditions}.
  \eitemize
\elemma
\bproof
  We first prove part (i).  Since $k\in\Ical$, it follows that $\|g^v_k\| > \epsilon$, which combined with~\eqref{eq.kkt-normal.dual} implies that $n_k \neq 0$, as claimed.  Now, in order to derive a contradiction, suppose that $0 = s_k = n_k + t_k$, which means that $-t_k = n_k \neq 0$.  From $\gv_k \neq 0$ and \eqref{eq.kkt-normal.dual}, it follows that $(\Hv_k + \lambdav_k I) n_k = -\gv_k \neq 0$, which gives
  \bequation\label{eq.zero_curvature}
    n_k^T (\Hv_k + \lambdav_k I) n_k = -n_k^T \gv_k = -n_k^T J_k^T c_k = -(J_k n_k)^T c_k = 0,
  \eequation
  where the last equality follows from $n_k = -t_k$ and $J_k t_k = 0$ (see~\eqref{eq.kkt_tangent.dual}). It now follows from~\eqref{eq.zero_curvature}, symmetry of $\Hv_k + \lambdav_k I$, and~\eqref{eq.kkt-normal.psd} that $0 = (\Hv_k + \lambdav_k I) n_k = -\gv_k$, which is a contradiction. This completes the proof of part (i).

  To prove part (ii), first observe that the conditions in \eqref{eq.tangential_step_conditions} are trivially satisfied if $t_k = 0$.  On the other hand, if Step~\ref{step.phase1.s_def} is reached with $t_k \neq 0$, then Step~\ref{step.phase1.tangential_reset} must have been reached, at which point it must have been determined that all of the conditions in \eqref{eq.tangential_step_conditions} held true (or else $t_k$ would have been reset to the zero vector).
\eproof

A key role in Algorithm~\ref{alg.phase1} is played by the sequence of trust funnel radii $\{\vmax_k\}$.  The next result establishes that it is a monotonically decreasing upper bound for the constraint violation, as previously claimed.

\blemma\label{lem.vmax}
  For all $k \in\Ical$, it follows that $v_k \leq \vmax_k$ and $0 < \vmax_{k+1} \leq \vmax_k$.
\elemma
\bproof
  The result holds trivially if $\Ical = \emptyset$.  Thus, let us assume that $\Ical \neq \emptyset$, which ensures that $0\in\Ical$.  Let us now use induction to prove the first inequality, as well as positivity of $\vmax_k$ for all $k \in \Ical$.  From the initialization in Algorithm~\ref{alg.phase1}, it follows that $v_0 \leq \vmax_0$ and $\vmax_0 > 0$.  Now, to complete the induction step, let us assume that $v_k \leq \vmax_k$ and $\vmax_k > 0$ for some $k\in\Ical$, then consider three cases.

  \textbf{Case 1: $k\in\Scal^f$.}  When $k \in \Scal^f$, let us consider the two possibilities based on the procedure for setting $\vmax_{k+1}$ stated in \eqref{eq.vmax_update_f}.  If \eqref{eq.vmax_update_f} sets $\vmax_{k+1} = v_{k+1} + \kappavtwo (\vmax_k - v_{k+1})$, then the fact that $k \in \Scal^f\subseteq\Fcal$, \eqref{eq.f-iter_vmax_rule}, and Lemma~\ref{lem.nonzero_steps}(i) imply that
  \bequationn
    \vmax_{k+1} = v_{k+1} + \kappavtwo(\vmax_k - v_{k+1}) \geq v_{k+1} + \kappavtwo \kapparho \|s_k\|^3 > v_{k+1} \geq 0.
  \eequationn
  On the other hand, if \eqref{eq.vmax_update_f} sets $\vmax_{k+1} = \max \{\kappavone \vmax_k, \vmax_k - \kapparho \|s_k\|^3\}$, then using the induction hypothesis, the fact that $k \in \Scal^f \subseteq \Fcal$, and \eqref{eq.f-iter_vmax_rule}, it follows that
  \bequationn
    \vmax_{k+1} \geq \kappavone \vmax_k > 0\ \ \text{and} \ \ \vmax_{k+1} \geq \vmax_k - \kapparho \|s_k\|^3 \geq v_{k+1} \geq 0.
  \eequationn
  This case is complete since, in each scenario, $\vmax_{k+1} \geq v_{k+1}$ and $\vmax_{k+1} > 0$.

  \textbf{Case 2: $k\in\Scal^v$.}  When $k \in \Scal^v$, let us consider the two possibilities based on the procedure for setting $\vmax_{k+1}$ stated in \eqref{eq.vmax_update_v}.  If \eqref{eq.vmax_update_v} sets $\vmax_{k+1} = v_{k+1} + \kappavtwo (\vmax_k - v_{k+1})$, then it follows from the induction hypothesis and the fact that $\rho_k^v \geq \kapparho$ for $k \in \Scal^v$ (which, in particular, implies that $v_{k+1} < v_k$ for $k \in \Scal^v$) that
  \bequationn
    \vmax_{k+1} = v_{k+1} + \kappavtwo(\vmax_k - v_{k+1}) \geq v_{k+1} + \kappavtwo(v_k - v_{k+1}) > v_{k+1} \geq 0.
  \eequationn
  On the other hand, if~\eqref{eq.vmax_update_v} sets $\vmax_{k+1} = \max \{\kappavone \vmax_k, v_{k+1} + \kappavtwo (v_k - v_{k+1})\}$, then the induction hypothesis and the fact that $v_{k+1} < v_k$ for $k \in \Scal^v$ implies that
  \bequationn
    \vmax_{k+1} \geq \kappavone \vmax_k > 0\ \ \text{and} \ \ \vmax_{k+1} \geq v_{k+1} + \kappavtwo(v_k - v_{k+1}) > v_{k+1} \geq 0.
  \eequationn
  This case is complete since, in each scenario, $\vmax_{k+1} \geq v_{k+1}$ and $\vmax_{k+1} > 0$.

  \textbf{Case 3: $k\notin\Scal^f\cup\Scal^v$.}  When $k \notin \Scal^f \cup \Scal^v$, it follows that $k\in\Ccal\cup\Ecal^v$, which may be combined with the induction hypothesis and the updating procedures for $x_k$ and $\vmax_k$ in Algorithms~\ref{alg.Fiter} and~\ref{alg.viter} to deduce that $0 < \vmax_k = \vmax_{k+1}$ and $v_{k+1} = v_k \leq \vmax_k = \vmax_{k+1}$.
  
  Combining the conclusions of the three cases above, it follows by induction that the first inequality of the lemma holds true and $\vmax_k > 0$ for all $k \in \Ical$.

  Let us now prove that $\vmax_{k+1} \leq \vmax_k$ for all $k \in \Ical$, again by considering three cases.  First, if $k \in \Scal^f$, then $\vmax_{k+1}$ is set using~\eqref{eq.vmax_update_f} such that
  \bequationn
    \vmax_{k+1} \leq \max\{\kappavone \vmax_k,\vmax_k - \kapparho \|s_k\|^3\} < \vmax_k,
  \eequationn
  where the strict inequality follows by $\kappavone \in (0,1)$ and Lemma~\ref{lem.nonzero_steps}(i).  Second, if $k \in \Scal^v$, then $v_{k+1} < v_k \leq \vmax_k$, where we have used the proved fact that $v_k \leq \vmax_k$; thus, $\vmax_k - v_{k+1} > 0$.  Then, since $\vmax_{k+1}$ is set using~\eqref{eq.vmax_update_v}, it follows that
  \bequationn
    \vmax_k - \vmax_{k+1} \geq \vmax_k - v_{k+1} - \kappavtwo (\vmax_k - v_{k+1}) = (1 - \kappavtwo) (\vmax_k - v_{k+1}) > 0.
  \eequationn
  Third, if $k \notin \Scal^f\cup\Scal^v$, then, by construction in Algorithms~\ref{alg.Fiter} and~\ref{alg.viter}, it follows that $\vmax_{k+1} = \vmax_k$.  This completes the proof.
\eproof

Our next lemma gives a lower bound for the decrease in the trust funnel radius as a result of a successful iteration.

\blemma\label{lem.vmax_decrease}
  If $k \in \Scal$, then $\vmax_k - \vmax_{k+1} \geq \kapparho(1-\kappavtwo)\|s_k\|^3$.
\elemma
\bproof
  If $k \in \Scal^f$, then $\vmax_{k+1}$ is set using~\eqref{eq.vmax_update_f}.  In this case,
  \bequalin
    \vmax_k - \vmax_{k+1}
      &\geq \vmax_k - v_{k+1} - \kappavtwo (\vmax_k - v_{k+1}) \\
      &=    (1-\kappavtwo) (\vmax_k - v_{k+1}) \geq \kapparho (1-\kappavtwo) \|s_k\|^3,
  \eequalin
  where the last inequality follows from~\eqref{eq.f-iter_vmax_rule} (since $k \in \Scal^f\subseteq\Fcal$).  If $k \in \Scal^v$, then $\vmax_{k+1}$ is set using~\eqref{eq.vmax_update_v}.  In this case, by Lemma \ref{lem.vmax} and the fact that $\rho_k^v \geq \kapparho$ for $k \in \Scal^v$,
  \bequalin
    \vmax_k - \vmax_{k+1} 
      &\geq \vmax_k - v_{k+1} - \kappavtwo (\vmax_k - v_{k+1}) \\
      &=    (1-\kappavtwo) (\vmax_k - v_{k+1}) \geq (1-\kappavtwo) (v_k - v_{k+1}) \geq \kapparho(1-\kappavtwo)\|s_k\|^3,
  \eequalin
  which completes the proof.
\eproof

Subsequently in our analysis, it will be convenient to consider an alternative formulation of problem~\eqref{prob.tangential_step} that arises from an orthogonal decomposition of the normal step $n_k$ into its projection onto the range space of $J_k^T$, call it $\NR_k$, and its projection onto the null space of $J_k$, call it $\NN_k$. Specifically, considering
\bequation\label{prob.tangential_step-rewritten2}
  \tNull_k \in \arg \min_{\tNull \in \R{N}} \mf_k (\NR_k+\tNull)\ \ \st \ \ J_k\tNull = 0\ \ \text{and}\ \ \|\tNull\| \leq \sqrt{(\deltas_k)^2 - \|\NR_k\|^2},
\eequation
we can recover the solution of \eqref{prob.tangential_step} as $t_k \gets \tNull_k - \NN_k$.  Similarly, for any $\lambdaf_k \in [0,\infty)$ that is strictly greater than the left-most eigenvalue of $Z_k^T H_k Z_k$, let us define
\bequation\label{prob.Q_f-rewritten}
  \bar{\Qcal}^f_k(\lambdaf_k) : \ \ \min_{\tNull \in \R{N}}\ (g_k + H_k \NR_k)^T \tNull + \thalf ({\tNull})^T(H_k+\lambdaf_k I) \tNull\ \ \st \ \ J_k \tNull = 0.
\eequation 
In the next lemma, we formally establish the equivalence between problems~\eqref{prob.tangential_step-rewritten2} and~\eqref{prob.tangential_step}, as well as between problems~\eqref{prob.Q_f-rewritten} and~\eqref{prob.Q_f}.

\blemma\label{lem.equivalence}
  For all $k \in \Ical$, the following problem equivalences hold:
  \bitemize
    \item[(i)]  if $\|n_k\| \leq \deltas_k$, then problems~\eqref{prob.tangential_step-rewritten2} and~\eqref{prob.tangential_step} are equivalent in that $(\tNull_k,\lambda_k^N)$ is part of a primal-dual solution of problem~\eqref{prob.tangential_step-rewritten2} if and only if $(t_k,\lambdaf_k) = (\tNull_k-\NN_k,\lambda_k^N)$ is part of a primal-dual solution of problem~\eqref{prob.tangential_step}; and
    \item[(ii)] if $Z_k^TH_kZ_k + \lambdaf_k I \succ 0$, then problems~\eqref{prob.Q_f-rewritten} and~\eqref{prob.Q_f} are equivalent in that $\tNull_k$ solves problem~\eqref{prob.Q_f-rewritten} if and only if $t_k = \tNull_k - \NN_k$ solves problem~\eqref{prob.Q_f}.
  \eitemize  
\elemma 
\bproof
  To prove part~(i), first note that $\|n_k\| \leq \deltas_k$ ensures that problems~\eqref{prob.tangential_step-rewritten2} and~\eqref{prob.tangential_step} are feasible.  Then, by $J_k\tNull = 0$ in \eqref{prob.tangential_step-rewritten2}, the vector $\NR_k \in \Range(J_k^T)$ is orthogonal with any feasible solution of \eqref{prob.tangential_step-rewritten2}, meaning that the trust region constraint in \eqref{prob.tangential_step-rewritten2} is equivalent to $\|\NR_k + \tNull\| \leq \deltas_k$.  Thus, as \eqref{eq.kkt_tangent} are the optimality conditions of \eqref{prob.tangential_step}, the optimality conditions of problem~\eqref{prob.tangential_step-rewritten2} (with this modified trust region constraint) are that there exists $(\tNull_k, y_k^N,\lambda_k^N) \in \R{N} \times \R{M} \times \R{}$ such that
  \bsubequations\label{eq.kkt_tangent-rewritten}
    \begin{align}
      \bbmatrix H_k + \lambda_k^N I & J_k^T \\ J_k & 0 \ebmatrix \bbmatrix \tNull_k \\ y_k^N \ebmatrix &= - \bbmatrix g_k + (H_k + \lambda_k^N I) \NR_k \\ 0 \ebmatrix, \label{eq.kkt_tangent-rewritten.dual} \\
      Z_k^TH_kZ_k + \lambda_k^N I &\succeq 0, \label{eq.kkt_tangent-rewritten.psd} \\ \text{and}\ \ 
      \lambda_k^N \perp (\deltas_k - \|\NR_k + \tNull_k\|) &\geq 0. \label{eq.kkt_tangent-rewritten.comp}
    \end{align}
  \esubequations
  From equivalence of the systems \eqref{eq.kkt_tangent-rewritten} and \eqref{eq.kkt_tangent}, it is clear that $(\tNull_k,y_k^N,\lambda_k^N)$ is a primal-dual solution of \eqref{prob.tangential_step-rewritten2} (with the modified trust region constraint) if and only if $(t_k,\yf_k,\lambdaf_k) = (\tNull_k-\NN_k,y_k^N,\lambda_k^N)$ is a primal-dual solution of \eqref{prob.tangential_step}.  This proves part (i).  Part (ii) follows in a similar manner from the orthogonal decomposition $n_k = \NN_k + \NR_k$ and the fact that $J_k\tNull=0$ in \eqref{prob.Q_f-rewritten} ensures that $\tNull_k \in \Null(J_k)$.
\eproof

The next lemma reveals important properties of the tangential step.  In particular, it shows that the procedure for performing a contraction of the trust region radius in an \Fiteration{} that results in a rejected step is well-defined.

\blemma\label{lem.wellposedness}
  If $k \in \Ccal^f$ and the condition in Step~\ref{step.lamcheck-Fiter} of Algorithm~\ref{alg.Fiter} tests true, then there exists $\lambdaf > \lambdaf_k$ such that $\sigmamin \leq \lambdaf / \|n_k + t(\lambdaf)\|$, where $t(\lambdaf)$ solves $\Qcal_k^f(\lambdaf)$.
\elemma
\bproof
  Since the condition in Step~\ref{step.lamcheck-Fiter} of Algorithm~\ref{alg.Fiter} is assumed to test true, it follows that $\lambdaf_k < \sigmamin \|s_k\|$.  Second, letting $t(\lambdaf)$ denote the solution of $\Qcal_k^f(\lambdaf)$, it follows by Lemma~\ref{lem.equivalence}(ii) that $\lim_{\lambdaf \to\infty} \|n_k + t(\lambdaf)\| = \|\NR_k\|$, meaning that $\lim_{\lambdaf \to\infty} \lambdaf / \|n_k + t(\lambdaf)\| = \infty$.
It follows from these observations and standard theory for trust region methods \cite[Chapter 7]{ConGT00a} that the result is true.
\eproof

The next lemma reveals properties of the normal step trust region radii along with some additional observations about the sequences $\{\deltamaxv_k\}$, $\{\lambdav_k\}$, and $\{\sigmav_k\}$.

\blemma \label{lem.deltav-changes}
  The following hold:
  \bitemize
    \item[(i)] if $k \in \Ccal^v$, then $0 < \deltav_{k+1} < \deltav_{k}$ and $\lambdav_{k+1} \geq \lambdav_k$;
    \item[(ii)] if $k \in \Ical$, then $\deltav_k \leq \deltamaxv_k \leq \deltamaxv_{k+1}$;
    \item[(iii)] if $k \in \Scal^v \cup \Ecal^v$, then $\deltav_{k+1} \geq \deltav_k$; and
    \item[(iv)] if $k \in \Fcal$, then $\deltav_{k+1} = \deltav_k$ and $\sigmav_{k+1} = \sigmav_k$. 
  \eitemize
\elemma
\bproof
  The proof of part (i) follows as that of~\cite[Lemma 3.4]{CurtRobiSama16}.  In particular, since the \text{V-}\contract\ procedure follows exactly that of \contract\ in \cite{CurtRobiSama16}, it follows that any call of \text{V-}\contract\ results in a contraction of the trust region radius for the normal subproblem and non-decrease of the corresponding dual variable.
  
  For part (ii), the result is trivial if $\Ical = \emptyset$.  Thus, let us assume that $\Ical \neq \emptyset$, which ensures that $0\in\Ical$.  We now first prove $\deltav_k \leq \deltamaxv_k$ for all $k \in \Ical$ using induction.  By the initialization procedure of Algorithm~\ref{alg.phase1}, it follows that $\deltav_0 \leq \deltamaxv_0$.  Hence, let us proceed by assuming that $\deltav_k \leq \deltamaxv_k$ for some $k \in \Ical$.  If $k \in \Scal^v$, then Step~\ref{step.delta_success-Viter} of Algorithm~\ref{alg.viter} shows that $\deltav_{k+1} \leq \deltamaxv_{k+1}$. If $k \in \Ecal^v$, then  Step~\ref{step.delta_expand-Viter} of Algorithm~\ref{alg.viter} gives $\deltav_{k+1} \leq \deltamaxv_{k+1}$. If $k \in \Ccal^v$, then part (i), Step~\ref{step.deltamaxC-Viter} of Algorithm~\ref{alg.viter}, and the induction hypothesis yield $\deltav_{k+1} < \deltav_k \leq \deltamaxv_k = \deltamaxv_{k+1}$. Lastly, if $k \in \Fcal$, then Step~\ref{step.phase1.f-iter-updates} of Algorithm~\ref{alg.phase1} and the inductive hypothesis give $\deltav_{k+1} = \deltav_k \leq \deltamaxv_k = \deltamaxv_{k+1}$. The induction step has now been completed since we have overall proved that $\deltav_{k+1} \leq \deltamaxv_{k+1}$, which means that we have proved the first inequality in part (ii). To prove $\deltamaxv_{k} \leq \deltamaxv_{k+1}$, consider two cases.  If $k \in \Scal^v$, then Step~\ref{step.deltamaxA-Viter} of Algorithm~\ref{alg.viter} gives $\deltamaxv_{k+1} \geq \deltamaxv_k$.  Otherwise, if $k \notin \Scal^v$, then according to Step~\ref{step.phase1.f-iter-updates} of Algorithm~\ref{alg.phase1} and Steps~\ref{step.deltamaxC-Viter} and \ref{step.deltamaxE-Viter} of Algorithm~\ref{alg.viter}, it follows that $\deltamaxv_{k+1} = \deltamaxv_k$.  Combining both cases, the proof of part (ii) is now complete.
  
  For part (iii), first observe from part (ii) and Step~\ref{step.delta_success-Viter} of Algorithm~\ref{alg.viter} that if $k \in \Scal^v$, then $\deltav_{k+1} = \min\{\deltamaxv_{k+1},\max\{\deltav_k, \gammae \|n_k\|\}\} \geq \deltav_k$.  On the other hand, if $k \in \Ecal^v$, then the conditions that must hold true for Step~\ref{step.expand-Viter} of Algorithm~\ref{alg.viter} to be reached ensure that $\lambdav_k > 0$, meaning that $\|n_k\| = \deltav_k$ (see~\eqref{eq.kkt-normal.comp}).  From this and the fact that the conditions in Step~\ref{step.expand-Viter} of Algorithm~\ref{alg.viter} must hold true, it follows that $\lambdav_k / \sigmav_k > \|n_k\| = \deltav_k$ and $\|n_k\| < \deltamaxv_k$.  Combining these observations with $\deltamaxv_{k+1} = \deltamaxv_k$ for $k \in \Ecal^v$ (see Step~\ref{step.deltamaxE-Viter} of Algorithm~\ref{alg.viter}) it follows from Step~\ref{step.delta_expand-Viter} of Algorithm~\ref{alg.viter} that $\deltav_{k+1} > \|n_k\| = \deltav_k$.

  Finally, part (iv) follows from Steps~\ref{step.phase1.f-iter-updates} and \ref{step.phase1.f-iter_sigmav} of Algorithm~\ref{alg.phase1}.
\eproof

The next result reveals similar properties for the other radii and $\{\lambdaf_k\}$.

\blemma \label{lem.deltaf-changes}
  The following hold:
  \bitemize
    \item[(i)] if $k \in \Ccal^f$, then $\deltaf_{k+1} < \deltaf_{k}$ and if, in addition, $(k+1) \in \Ical^t$, then $\lambdaf_{k+1} \geq \lambdaf_k$;
    \item[(ii)] if $k \in \Scal^f$, then $\deltaf_{k+1} \geq \deltaf_k$ and $\deltas_{k+1} \geq \deltas_k$; and
    \item[(iii)] if $k \in \Vcal$, then $\deltaf_{k+1} = \deltaf_k$.
  \eitemize
\elemma
\bproof
  For part (i), notice that $\deltaf_{k+1}$ is set in Step~\ref{step.delta_contract-Fiter} of Algorithm~\ref{alg.Fiter} and that $(x_{k+1},\deltav_{k+1}) \gets (x_k,\deltav_k)$ and $n_{k+1} = n_k$ for all $k \in \Ccal^f$.  Let us proceed by considering two cases depending on the condition in Step~\ref{step.lamcheck-Fiter} of Algorithm~\ref{alg.Fiter}.
 
  \textbf{Case 1:}  $\lambdaf_k < \sigmamin\|s_k\|$.  In this case, $\deltaf_{k+1}$ is set in Step~\ref{step.deltamagic-Fiter} of Algorithm~\ref{alg.Fiter}, which from Step~\ref{step.linsysmagic-Fiter} of Algorithm~\ref{alg.Fiter} and Lemma~\ref{lem.wellposedness} implies that $\lambdaf > \lambdaf_k$.  Combining this with Lemma~\ref{lem.equivalence} and standard theory for trust region methods leads to the fact that the solution $\tNull(\lambdaf)$ of $\bar{\Qcal}^f_k(\lambdaf)$ satisfies $\|\tNull(\lambdaf)\| < \|\tNull_k\|$.  Thus, $\deltaf_{k+1} = \|n_k + t(\lambdaf)\| = \|\NR_k + \tNull(\lambdaf)\| < \|\NR_k+\tNull_k\| = \|s_k\| \leq \deltaf_k$, where the last inequality comes from~\eqref{def.delta_s}.  If, in addition, $(k+1) \in \Ical^t$ so that a nonzero tangential step is computed and not reset to zero, it follows that $\lambdaf_{k+1} = \lambdaf$.  This establishes the last conclusion of part~(i) for this case since it has already been shown above that $\lambdaf > \lambdaf_k$.
 
  \textbf{Case 2:} $\lambdaf_k \geq \sigmamin\|s_k\|$.  In this case, $\deltaf_{k+1}$ is set in Step~\ref{step.deltanormal-Fiter} of Algorithm~\ref{alg.Fiter} and, from~\eqref{def.delta_s} and $\gammac \in (0,1)$, it follows that $\deltaf_{k+1} = \gammac \|s_k\| \leq \gammac \deltaf_k < \deltaf_k$.  Consequently, from Step~\ref{step.phase1.f-iter-updates} of Algorithm~\ref{alg.phase1} and~\eqref{def.delta_s}, one finds that $\deltas_{k+1} \leq \deltas_k$.  It then follows from Lemma~\ref{lem.equivalence} and standard trust region theory that if $(k+1) \in \Ical^t$, then $\lambdaf_{k+1} \geq \lambdaf_k$.

  To prove part (ii), notice that for $k \in \Scal^f$ it follows by Step~\ref{step.delta_success-Fiter} of Algorithm~\ref{alg.Fiter} that $\deltaf_{k+1} = \max\{\deltaf_k,\gammae\|s_k\|\}$, so $\deltaf_{k+1} \geq \deltaf_k$. From this, Step~\ref{step.phase1.f-iter-updates} of Algorithm~\ref{alg.phase1}, and~\eqref{def.delta_s} it follows that $\deltas_{k+1} \geq \deltas_k$.  These conclusions complete the proof of part (ii).

  Finally, part (iii) follows from Step~\ref{step.phase1.v-iter-updates} of Algorithm~\ref{alg.phase1}.
\eproof

Next, we show that after a \Viteration{} with either a contraction or an expansion of the trust region radius, the subsequent iteration cannot result in an expansion.

\blemma\label{lem.not-expansion}
  If $k \in \Ccal^v \cup \Ecal^v$, then $(k+1) \in \Fcal \cup \Scal^v \cup \Ccal^v$.
\elemma
\bproof
  If $(k + 1) \in \Fcal$, then there is nothing left to prove.  Otherwise, if $(k + 1) \in \Vcal$, then the proof follows using the same logic as for~\cite[Lemma 3.7]{CurtRobiSama16}, which shows that one of three cases holds: $(i)$ $k \in \Ccal^v$, which yields $\lambdav_{k+1} \leq \sigmav_{k+1}\|n_{k+1}\|$, so $(k+1) \notin \Ecal^v$; $(ii)$ $k \in \Ecal^v$ and $\deltamaxv_k \geq \lambdav_k/\sigmav_k$, which also yields $\lambdav_{k+1} \leq \sigmav_{k+1}\|n_{k+1}\|$, so $(k+1) \notin \Ecal^v$; or $(iii)$ $k \in \Ecal^v$ and $\deltamaxv_k < \lambdav_k/\sigmav_k$, which implies $(k+1) \in \Scal^v \cup \Ccal^v$, so $(k+1) \notin \Ecal^v$.
\eproof

Our goal now is to expand upon the conclusions of Lemma \ref{lem.not-expansion}.  To do this, it will be convenient to define the first index in a given index set following an earlier index~$\kbar \in \Ical$ in that index set (or the initial index 0).  In particular, let us define
\begin{align*}
  k_{\Scal}(\kbar) := \min\{k \in \Scal : k > \kbar\}\ \ 
  \text{and}\ \ 
  k_{\Scal\cup\Vcal}(\kbar) := \min\{k \in \Scal \cup \Vcal : k > \kbar\}
\end{align*}
along with the associated sets
\begin{align*}
  \Ical_{\Scal}(\kbar) := \{k \in \Ical :\kbar < k < k_{\Scal}(\kbar)\}\ \ \text{and}\ \ 
  \Ical_{\Scal\cup\Vcal}(\kbar) := \{k \in \Ical :\kbar < k < k_{\Scal\cup\Vcal}(\kbar)\}.
\end{align*}

The following lemma shows one important property related to these quantities.

\blemma\label{lem.1-expansion}
  For all $\kbar \in \Scal \cup \{0\}$, it follows that $| \Ecal^v \cap \Ical_\Scal(\kbar) | \leq 1$.
\elemma
\bproof
  In order to derive a contradiction, suppose that there exists $\kbar \in \Scal \cup \{0\}$ such that $|\Ecal^v \cap \Ical_\Scal(\kbar)| > 1$, which means that one can choose $\kSone$ and $\kSthree$ as the first two distinct indices in $\Ecal^v \cap \Ical_\Scal(\kbar)$; in particular,
  \bequationn
    \{\kSone,\kSthree\}\subseteq \Ecal^v \cap \Ical_\Scal(\kbar)\ \ \text{and} \ \ \kbar < \kSone < \kSthree < k_\Scal(\kbar).
  \eequationn
  By Lemma~\ref{lem.not-expansion} and the fact that $\kSone \in \Ecal^v$, it follows that $\{\kSone+1, \dots, \kSthree-1\} \neq \emptyset$.  Let us proceed by considering two cases, deriving a contradiction in each case.
  
  \textbf{Case 1:} $\Vcal \cap \{\kSone+1, \dots, \kSthree-1\} = \emptyset$.  In this case, by the definitions of $\kSone$, $\kSthree$, and $\Ical_\Scal(\kbar)$, it follows that $\{\kSone+1, \dots\kSthree-1\} \subseteq \Ccal^f$.  Then, since $\deltav_{k+1} = \deltav_k$ and $\sigmav_{k+1} = \sigmav_k$ for all $k \in \Ccal^f \subseteq \Fcal$, it follows that $\deltav_{\kSthree} = \deltav_{\kSone+1}$ and $\sigmav_{\kSthree} = \sigmav_{\kSone+1}$.  In particular, using the fact that $\deltav_{\kSthree} = \deltav_{\kSone+1}$, it follows along with the fact that $x_{k+1} = x_k$ for all $k \notin \Scal$ that $\|n_{\kSthree}\| = \|n_{\kSone+1}\|$ and $\lambdav_{\kSthree} = \lambdav_{\kSone+1}$.  Now, since $(\kSone+1) \in \Ccal^f$, it follows with Step~\ref{step.phase1.check_type} of Algorithm~\ref{alg.phase1} and \eqref{eq.f-iter_sigma_v} that
  \bequationn
    \lambdav_{\kSthree}/\|n_{\kSthree}\| = \lambdav_{\kSone+1}/\|n_{\kSone+1}\| \leq \sigmav_{\kSone+1} = \sigmav_{\kSthree},
  \eequationn
  which implies that $\kSthree \notin \Ecal^v$, a contradiction.
%
%
  
  \textbf{Case 2:} $\Vcal \cap \{\kSone+1, \dots, \kSthree-1\} \neq \emptyset$.  In this case, by the definitions of $\kSone$, $\kSthree$, and $\Ical_\Scal(\kbar)$, it follows that $\{\kSone+1, \dots, \kSthree-1\} \subseteq \Ccal^f \cup \Ccal^v$.  In addition, by the condition of this case, it also follows that there exists a greatest index $\kStwo\in\Ccal^v \cap \{\kSone+1, \dots, \kSthree-1\}$.  In particular, for the index $\kStwo \in \Ccal^v$, it follows that $\kSone+1 \leq \kStwo \leq \kSthree-1$ and $\{\kStwo+1,\dots,\kSthree-1\} \subseteq \Ccal^f$.  By $\kStwo \in \Ccal^v$ and Lemma~\ref{lem.not-expansion}, it follows that $\kStwo+1 \notin \Ecal^v$; hence, since $\kSthree \in \Ecal^v$, it follows that $\{\kStwo+1,\dots, \kSthree-1\} \neq \emptyset$.  We may now apply the same argument as for Case 1, but with $\kSone$ replaced by $\kStwo$, to arrive at the contradictory conclusion that $\kSthree\notin\Ecal^v$, completing the proof.
\eproof

The next lemma reveals lower bounds for the norms of the normal and full steps.

\blemma\label{lem.bound-on-normn}
  For all $k \in \Ical$, the following hold:
  \bitemize
    \item[(i)] $\|n_k\| \geq \min \left\{ \deltav_k, \|\gv_k\|/ \|\Hv_k\| \right\} > 0$ and
    \item[(ii)] $\|s_k\| \geq \kappantn \min \left\{ \deltav_k, \|\gv_k\| / \|\Hv_k\| \right\} > 0$.
  \eitemize
\elemma
\bproof
  The proof of part (i) follows as that for~\cite[Lemma 3.2]{CurtRobiSama16}.
Part (ii) follows from part (i) and \eqref{eq.s_vs_n}, the latter of which holds because of Lemma~\ref{lem.nonzero_steps}(ii).
\eproof

We now provide a lower bound for the decrease in the model of infeasibility.

\blemma \label{lem.vmodeldec-1}
  For all $k\in\Ical$, the quantities $n_k$, $\lambdav_k$, and $s_k$ satisfy
  \bsubequations
    \begin{align}
      v_k - \mv_k(n_k) &= \thalf n_k^T (\Hv_k + \lambdav_k I) n_k + \thalf \lambdav_k \|n_k\|^2 > 0,  \label{eq.vmodeldec} \\
      v_k - \mv_k(s_k) &\geq \kappavm (\thalf n_k^T (\Hv_k+\lambdav_k I)n_k + \thalf \lambdav_k \|n_k\|^2) > 0, \ \ \text{and} \label{eq.vmodeldec-s} \\
      v_k - \mv_k(s_k) &\geq \thalf \kappavm \|\gv_k\| \min \left\{\deltav_k, \|\gv_k\| / \|\Hv_k\| \right\} > 0.  \label{eq.vmodeldec-s-2}
    \end{align}
  \esubequations
\elemma
\bproof
  The proof of~\eqref{eq.vmodeldec} follows as for that of~\cite[Lemma 3.3]{CurtRobiSama16} and the fact that $\|\gv_k\| > \eps$, which holds since $k\in\Ical$.  The inequalities in~\eqref{eq.vmodeldec-s} follow from~\eqref{eq.vmodeldec} and \eqref{eq.mv_red_t}, the latter of which holds because of Lemma~\ref{lem.nonzero_steps}(ii). To prove~\eqref{eq.vmodeldec-s-2}, first observe from standard trust region theory (e.g., see~\cite[Theorem~6.3.1]{ConGT00a}) that
  \bequation\label{eq.v-modelDecrease}
v_k - \mv_k(n_k) \geq \thalf\|\gv_k\| \min \left\{\deltav_k, \|\gv_k\| / \|\Hv_k\| \right\} > 0.
  \eequation
  By combining~\eqref{eq.v-modelDecrease} and~\eqref{eq.mv_red_t} (which holds by Lemma~\ref{lem.nonzero_steps}(ii)), one obtains~\eqref{eq.vmodeldec-s-2}.
\eproof

The next lemma reveals that if the dual variable for the normal step trust region is beyond a certain threshold, then the trust region constraint must be active and the step will either be an \Fiteration\ or a successful \Viteration.  Consequently, this reveals an upper bound for the dual variable for any unsuccessful \Viteration.

\blemma\label{lem.largelambdav->acceptedstep}
  For all $k \in \Ical$, if the trial step $s_k$ and the dual variable $\lambdav_k$ satisfy
  \bequation\label{eq.largelambdav->acceptedstep}
    \lambdav_k \geq \frac{\kappadelta^2}{\kappavm}(2\gv_{Lip}+\thetafc+2\kapparho\|s_k\|),
  \eequation
  then $\|n_k\| = \deltav_k$ and $\rho_k^v \geq \kapparho$.
\elemma
\bproof
  For all $k \in \Ical$, it follows from the definition of $\mv_k$ and the Mean Value Theorem that there exists a point $\bar{x}_k \in \R{N}$ on the line segment $[x_k, x_k + s_k]$ such that
  \begin{align}
    \mv_k(s_k) - v(x_k + s_k)
      &= \(\gv_k - \gv(\bar{x}_k)\)^T s_k + \thalf s_k^T \Hv_k s_k \nonumber \\
      &\geq - \| \gv_k - \gv(\bar{x}_k) \| \|s_k\| - \thalf \|\Hv_k\| \|s_k\|^2. \label{eq.largelambdav->acceptedstep-eq2}
  \end{align}
  By \eqref{eq.largelambdav->acceptedstep} and \eqref{eq.kkt-normal.comp}, it follows that $\|n_k\| = \deltav_k$.  Combining this fact with~\eqref{eq.largelambdav->acceptedstep-eq2}, \eqref{eq.vmodeldec-s}, \eqref{eq.kkt-normal.psd}, Lemma~\ref{lem.derivative_bounds}, \eqref{eq.largelambdav->acceptedstep}, and the fact that $\|s_k\| \leq \deltas_k \leq \kappadelta \deltav_k = \kappadelta \|n_k\|$, one obtains
  \bequalin
    v_k - v(x_k + s_k) &= v_k - \mv_k(s_k) + \mv_k(s_k) - v(x_k + s_k) \\
    &\geq \thalf \kappavm \lambdav_k \|n_k\|^2 - \| \gv_k - \gv(\bar{x}_k) \| \|s_k\| - \thalf \|\Hv_k\| \|s_k\|^2 \\
    &\geq \thalf \kappavm \kappadelta^{-2} \lambdav_k \|s_k\|^2 - \| \gv_k - \gv(\bar{x}_k) \| \|s_k\| - \thalf \|\Hv_k\| \|s_k\|^2 \\
    &\geq (\thalf \kappavm \kappadelta^{-2} \lambdav_k - \gv_{Lip} - \thalf \thetafc) \|s_k\|^2 
    \geq \kapparho \|s_k\|^3,
  \eequalin
  which, by Steps~\ref{step.phase1.f-iter-updates} and \ref{step.phase1.v-iter-rho} in Algorithm~\ref{alg.phase1} and \eqref{def.rho_v}, completes the proof.
\eproof

Recall that our main goal in this section is to prove that $|\Ical| < \infty$.  Ultimately, this result is attained by deriving contradictions under the assumption that $|\Ical| = \infty$.  For example, if $|\Ical| = \infty$ and the iterations corresponding to all sufficiently large $k \in \Ical$ involve contractions of a trust region radius, then the following lemma helps to lead to contradictions in subsequent results.  In particular, it reveals that, under these conditions, a corresponding dual variable tends to infinity.

\blemma \label{lem.lambdav->infty}
  The following hold:
  \begin{itemize}
    \item[(i)] If $k \notin \Scal$ for all large $k \in \Ical$ and $|\Ccal^v| = \infty$, then $\{\deltav_k\} \to 0$ and $\{\lambdav_k\} \to \infty$.
    \item[(ii)] If $k \in \Ccal^f$ for all large $k \in \Ical$, then $\{\deltaf_k\} \to 0$ and $\{\lambdaf_k\} \to \infty$.
  \end{itemize}
\elemma
\bproof
  By Lemma~\ref{lem.deltav-changes}, Lemma~\ref{lem.1-expansion}, and the fact that $k\notin\Scal$ for all large $k \in \Ical$, the proof of part (i) follows as that of~\cite[Lemma 3.9]{CurtRobiSama16}.

  To prove part (ii), let us assume, without loss of generality, that $k \in \Ccal^f$ for all $k \in \Ical$.  It then follows that $k \in \Ical^t$ for all $k \in \Ical$, since otherwise it would follow that $t_k \gets 0$, which by \eqref{eq.f-iter_t_vs_s} means $k \in \Vcal$, a contradiction to $k \in \Ccal^f$.  Thus,
  \bequation\label{eq:CcapI}
    k \in \Ccal^f \cap \Ical^t\ \ \text{for all}\ \ k \in\Iscr.
  \eequation
  Next, we claim that the condition in Step~\ref{step.lamcheck-Fiter} of Algorithm~\ref{alg.Fiter} can hold true for at most one iteration.  If it never holds true, then there is nothing left to prove.  Otherwise, let $k_c \in \Ical$ be the first index for which the condition holds true.  The structure of Algorithm~\ref{alg.Fiter} (see Step~\ref{step.linsysmagic-Fiter}) and~\eqref{eq:CcapI} then ensure that $\lambdaf_{k_c+1} / \|s_{k_c+1}\| \geq \sigmamin$. From Lemma~\ref{lem.deltaf-changes}(i), one may conclude that $\{\lambdaf_k / \|s_k\|\}$ is nondecreasing.  From this, it follows that the condition in Step~\ref{step.lamcheck-Fiter} of Algorithm~\ref{alg.Fiter} will never be true for any $k > k_c$. Thus, we may now proceed, without loss of generality, under the assumption that the condition in Step~\ref{step.lamcheck-Fiter} of Algorithm~\ref{alg.Fiter} always tests false.  This means that $\deltaf_{k+1}$ is set in Step~\ref{step.deltanormal-Fiter} of Algorithm~\ref{alg.Fiter} for all $k\in\Ical$, yielding $\deltaf_{k+1} \gets \gammac \|s_k\| \leq \gammac \deltaf_k$, where the last inequality comes from~\eqref{def.delta_s}. Therefore, $\{\deltaf_k\} \to 0$ for all $k\in\Ical$, and consequently $\{\lambdaf_k\} \to \infty$.
\eproof

We now show that the sequences $\{\deltamaxv_k\}$ and $\{n_k\}$ are bounded above.

\blemma \label{lem.deltamax-bound}
  There exists a scalar $\deltamaxv_{\max} \in (0,\infty)$ such that $\deltamaxv_k = \deltamaxv_{\max}$ for all sufficiently large $k \in \Ical$. In addition, $|\Scal^v_{\deltamax}| < \infty$ and there exists a scalar $n_{\max} \in (0,\infty)$ such that $\|n_k\| \leq n_{\max}$ for all $k \in \Ical$.
\elemma
\bproof
First, in order to derive a contradiction, assume that there is no $\deltamaxv_{\max}$ such that $\deltamaxv_k = \deltamaxv_{\max}$ for all sufficiently large $k \in \Ical$.  This, in turn, means that Step~\ref{step.deltamaxA-Viter} of Algorithm~\ref{alg.viter} is reached infinitely often, meaning that $|\Scal^v| = \infty$. For all $k \in \Scal^v \subseteq \Scal$, it follows from Lemma~\ref{lem.vmax_decrease} that $\vmax_k - \vmax_{k+1} \geq \kapparho(1-\kappavtwo)\|s_k\|^3$. Now, using the monotonicity of $\{\vmax_k\}$ and the fact that $\vmax_k \geq 0$ (see Lemma~\ref{lem.vmax}), one may conclude that $\{\vmax_k\}$ converges; therefore $\{s_k\}_{k\in\Scal^v} \to 0$. From this fact, Lemma~\ref{lem.nonzero_steps}(ii), and~\eqref{eq.s_vs_n} it follows that $\{n_k\}_{k\in\Scal^v} \to 0$. Thus, there exists an iteration index $k^v_{\deltamax}$ such that for all $k \in \Scal^v$ with $k \geq k^v_{\deltamax}$, one finds $\gammae \|n_k\| < \deltamaxv_0 \leq \deltamaxv_k$, where the last inequality follows from Lemma~\ref{lem.deltav-changes}(ii). From this and Steps~\ref{step.deltamaxA-Viter}, \ref{step.deltamaxC-Viter}, and \ref{step.deltamaxE-Viter} of Algorithm~\ref{alg.viter}, it follows that $\deltamaxv_{k+1} \gets \deltamaxv_k$ for all $k \geq k^v_{\deltamax}$, a contradiction.  The proof of the second part of the lemma follows as in that for~\cite[Lemma 3.11]{CurtRobiSama16}.
\eproof

In the next lemma, a uniform lower bound on $\{\deltav_k\}$ is provided.

\blemma \label{lem.lowerboundondeltav}
  There exists a scalar $\deltav_{\min} \in (0,\infty)$ such that $\deltav_k \geq \deltav_{\min}$ for all $k \in \Ical$.
\elemma
\bproof
  If $|\Ccal^v| < \infty$, then the result follows from Lemma~\ref{lem.deltav-changes}(iii)--(iv).  Thus, let us proceed under the assumption that $|\Ccal^v| = \infty$.  As in the beginning of the proof of Lemma~\ref{lem.largelambdav->acceptedstep}, it follows that~\eqref{eq.largelambdav->acceptedstep-eq2} holds.  Then, using~\eqref{eq.largelambdav->acceptedstep-eq2}, \eqref{eq.vmodeldec-s-2}, Lemma~\ref{lem.derivative_bounds}(i), Lemma~\ref{lem.derivative_bounds}(iii), $\|\gv_k\| >\epsilon$ for $k\in\Ical$, and $\|s_k\| \leq \deltas_k \leq \kappadelta \deltav_k$, it follows that
  \bequalin
    &\phantom{=l} v_k - v(x_k + s_k) \\
    &= v_k - \mv_k(s_k) + \mv_k(s_k) - v(x_k + s_k) \\
    &\geq \thalf \kappavm \|\gv_k\| \min \left\{\deltav_k, \|\gv_k\| / \|\Hv_k\| \right\} - \| \gv_k - \gv(\bar{x}_k) \| \|s_k\| - \thalf \|\Hv_k\| \|s_k\|^2 \\
    &\geq \thalf \kappavm \epsilon \min \left\{\deltav_k, \epsilon / \thetafc \right\} - (\gv_{Lip} + \thalf\thetafc) \|s_k\|^2 \\
    &\geq \thalf \kappavm \epsilon \min \left\{\deltav_k, \epsilon / \thetafc \right\} - (\gv_{Lip} + \thalf\thetafc) \kappadelta^2 (\deltav_k)^2.
  \eequalin
  Considering these inequalities and $\|s_k\| \leq \deltas_k \leq \kappadelta \deltav_k$, it must hold that $\rho_k^v \geq \kapparho$ for any $k \in \Ical$ as long as $\deltav_k \in (0, \epsilon / \thetafc]$ is sufficiently small such that
  \bequationn
    \thalf \kappavm \epsilon \deltav_k - (\gv_{Lip} + \thalf \thetafc) \kappadelta^2 (\deltav_k)^2 \geq \kapparho \kappadelta^3 (\deltav_k)^3 \geq \kapparho \|s_k\|^3.
  \eequationn
  This fact implies the existence of a positive threshold $\deltav_{thresh} \in (0, \epsilon / \thetafc]$ such that, for any $k \in \Ical$ with $\deltav_k \in (0, \deltav_{thresh})$, one finds $\rho_k^v \geq \kapparho$.  Along with the fact that $\rho_k^v < \kapparho$ if and only if $k\in\Ccal^v$ (see Step~\ref{step.accept-Viter}, \ref{step.contract-Viter}, and \ref{step.expand-Viter} of Algorithm~\ref{alg.viter} and Step~\ref{step.phase1.f-iter-updates} of Algorithm~\ref{alg.phase1}), it follows that
  \begin{equation} \label{eq.key-Cv}
    \deltav_k \geq \deltav_{thresh} \ \text{for all $k\in\Ccal^v$.}
  \end{equation}
  Since the normal step subproblem trust region radius is only decreased when $k \in \Ccal^v$, we will complete the proof by showing a lower bound on $\deltav_{k+1}$ when $k\in\Ccal^v$.

  Suppose that $k \in \Ccal^v$.  If Step~\ref{step.regularized1-Viter} of Algorithm~\ref{alg.viter} is reached, then 
  \bequationn
    \deltav_{k+1} \gets \|n(\lambdav)\| \geq \frac{\lambdav}{\sigmamax} = \frac{\lambdav_k + (\sigmamin \|\gv_k\|)^{1/2}}{\sigmamax} \geq \frac{(\sigmamin \|\gv_k\|)^{1/2}}{\sigmamax} \geq \frac{(\sigmamin \epsilon)^{1/2}}{\sigmamax},
  \eequationn
  where the last inequality follows since $k \in \Ical$ means $\|\gv_k\| \geq \epsilon$.  If Step~\ref{step.deltamagic-Viter} is reached, then the algorithm chooses $\lambdav \in (\lambdav_k,\hat{\lambda}^v)$ to find $n(\lambdav)$ that solves $\mathcal{Q}_k^v(\lambdav)$ such that $\sigmamin \leq \lambdav/\|n(\lambdav)\| \leq \sigmamax$. For this case and the cases when Step~\ref{step.delta2-Viter} or \ref{step.deltanormal-Viter} is reached, the existence of $\deltav_{\min} \in (0,\infty)$ such that $\deltav_{k+1} \geq \deltav_{\min}$ for all $k \in \Ccal^v$ follows in the same manner as in the proof of \cite[Lemma 3.12]{CurtRobiSama16}.  Combining these facts with~\eqref{eq.key-Cv} and Lemma~\ref{lem.deltav-changes}(iii)--(iv), the proof is complete.
\eproof

The next result shows that there are finitely many successful iterations.

\blemma\label{lem.vmax-->0}
  The following hold: $|\Scal^v| < \infty$ and $|\Scal^f| < \infty$.
\elemma
\bproof
  Lemma~\ref{lem.lowerboundondeltav}, $\|g^v_k\| > \epsilon$ for all $k\in\Ical$, Lemma~\ref{lem.bound-on-normn}(i), and Lemma~\ref{lem.derivative_bounds}(i) imply the existence of $n_{\min} \in (0,\infty)$ such that $\|n_k\| \geq n_{\min}$ for all $k\in\Ical$, i.e.,
  \bequation\label{summarize}
    \|g^v_k\| > \epsilon \ \ \text{and} \ \ \|n_k\| \geq n_{\min} > 0 \ \ \text{for all $k\in\Ical$.}
  \eequation
  In order to reach a contradiction to the first desired conclusion, suppose that $|\Scal^v| = \infty$.  For any $k \in \Scal^v$, it follows from Lemma~\ref{lem.vmax_decrease}, Lemma~\ref{lem.nonzero_steps}(ii), and \eqref{eq.s_vs_n} that
  \bequation\label{vmak-change}
    \vmax_k - \vmax_{k+1} \geq \kapparho (1-\kappavtwo) \|s_k\|^3\geq \kapparho (1-\kappavtwo) \kappantn^3 \|n_k\|^3.
  \eequation
  By Lemma~\ref{lem.vmax}, $0 < \vmax_{k+1} \leq \vmax_k$ for all $k \in \Ical$, meaning that $\{\vmax_k - \vmax_{k+1}\} \to 0$, which together with~\eqref{vmak-change} shows that $\{\|n_k\|\}_{k \in \Scal^v} \to 0$, contradicting~\eqref{summarize}.  This proves that $|\Scal^v| < \infty$.  Now, in order to reach a contradiction to the second desired conclusion, suppose that $|\Scal^f| = \infty$.  Since $|\Scal^v| < \infty$, we can assume without loss of generality that $\Scal = \Scal^f$.  This means that the sequence $\{f_k\}$ is monotonically nonincreasing.  Combining this with the fact that $\{f_k\}$ is bounded below under Assumptions~\ref{ass.differentiability} and~\ref{ass.compactness}, it follows that $\{f_k\} \to f_{low}$ for some $f_{low} \in (-\infty,\infty)$ and $\{f_k - f_{k+1}\} \to 0$.  Using these facts, the inequality $\rho^f_k \geq \kapparho$ for all $k\in\Scal^f$, and $|\Scal^f| = \infty$, it follows that $\{\kapparho \|s_k\|^3\}_{k \in \Scal^f} \leq \{f_k - f_{k+1}\}_{k \in \Scal^f} \to 0$, which gives $\{\|s_k\|\}_{k \in \Scal^f} \to 0$.  This, in turn, implies that $\{\|n_k\|\}_{k \in \Scal^f} \to 0$ because of Lemma~\ref{lem.nonzero_steps}(ii) and~\eqref{eq.s_vs_n}, which contradicts~\eqref{summarize}.  Hence, $|\Scal^f| < \infty$.
\eproof

We are now prepared to prove that Algorithm~\ref{alg.phase1} terminates finitely.

\btheorem \label{thm.finitetermination}
  Algorithm~\ref{alg.phase1} terminates finitely, i.e., $|\Ical| < \infty$.
\etheorem
\bproof
  Suppose by contradiction that $|\Ical| = \infty$.  
  Let us consider two cases.
  
  \textbf{Case 1:} $|\Vcal| = \infty$.  Since $|\Scal| < \infty$, it follows that $|\Vcal \setminus \Scal^v| = |\Ccal^v\cup\Ecal^v| = \infty$, which along with Lemma~\ref{lem.1-expansion} implies that $|\Ecal^v| < \infty$ while $|\Ccal^v| = \infty$.  It now follows from Lemma~\ref{lem.lambdav->infty}(i) that $\{\deltav_k\} \to 0$,
which contradicts Lemma~\ref{lem.lowerboundondeltav}.

  \textbf{Case 2:} $|\Vcal| < \infty$.  For this case, we may assume without loss of generality that $\Fcal = \Ical$.  This implies with Lemma~\ref{lem.nonzero_steps}(i) that $\deltav_k = \deltav_0$ and $n_k = n_0 \neq 0$ for all $k \in \Ical$.  It also implies from Step~\ref{step.phase1.check_type} of Algorithm~\ref{alg.phase1} that~\eqref{eq.f-iteration-conditions} holds for all $k \in \Ical$; in particular, from~\eqref{eq.f-iter_t_vs_s} it means that $t_k \neq 0$ for all $k \in \Ical$.  Now, from $|\Vcal| < \infty$, $|\Scal| < \infty$, and Lemma~\ref{lem.lambdav->infty}(ii), it follows that $\{\deltaf_k\} \to 0$, which by~\eqref{def.delta_s} yields $\{\deltas_k\} \to 0$.  It then follows from Step~\ref{step.phase1.tangential_conditions} of Algorithm~\ref{alg.phase1} and $\Fcal = \Ical$ that $\{n_k\} \to 0$, which contradicts our previous conclusion that $n_k = n_0 \neq 0$ for all $k \in \Ical$.
\eproof

\subsection{Complexity analysis for phase 1}\label{sec.complexity-feasible}

Our goal in this subsection is to prove an upper bound on the total number of iterations required until phase 1 terminates, i.e., until the algorithm reaches $k \in \N{}$ such that $\|\gv_k\| \leq \epsilon$.  To prove such a bound, we require the following additional assumption.

\bassumption \label{ass.M-lip}
  The Hessian functions $\Hv(x) := \nabla^2 v(x)$ and $\nabla^2 f(x)$ are Lipschitz continuous with constants $\Hv_{Lip} \in (0,\infty)$ and $H_{Lip} \in (0,\infty)$, respectively, on a path defined by the sequence of iterates and trial steps computed in the algorithm.
\eassumption

Our first result in this subsection can be seen as a similar conclusion to that given by Lemma~\ref{lem.largelambdav->acceptedstep}, but with this additional assumption in hand.

\blemma \label{lem.lambdak-lowerbound-1}
  For all $k \in \Ical$, if the trial step $s_k$ and dual variable $\lambdav_k$ satisfy
  \bequation \label{eq.largeenoughlambda-1}
    \lambdav_k \geq \kappadelta^2\kappavm^{-1}(\Hv_{Lip} + 2\kapparho)\|s_k\|,
\eequation
  then $\|n_k\| = \deltav_k$ and $\rho_k^v \geq \kapparho$.
\elemma
\bproof
  For all $k \in \Ical$, there exists $\overline{x}_k$ on the line segment $[x_k, x_k + s_k]$ such that
  \bequation
    \mv_k(s_k) - v(x_k + s_k) = \thalf s_k^T \(\Hv_k - \Hv(\overline{x}_k)\) s_k \geq -\thalf \Hv_{Lip} \|s_k\|^3.
  \eequation
  From this, \eqref{eq.vmodeldec-s}, and \eqref{eq.kkt-normal.psd}, one deduces that
  \bequalin
    v(x_k) - v(x_k + s_k) &= v(x_k) - \mv_k(s_k) + \mv_k(s_k) - v(x_k + s_k)\\
      &\geq \thalf\kappavm \lambdav_k \|n_k\|^2 - \thalf \Hv_{Lip} \|s_k\|^3.
  \eequalin
  From Lemma~\ref{lem.nonzero_steps}(i), \eqref{eq.largeenoughlambda-1}, and~\eqref{eq.kkt-normal.comp}, it follows that $\|n_k\| = \deltav_k$, which along with~\eqref{def.delta_s} means that $\|s_k\| \leq \deltas_k \leq \kappadelta \deltav_k = \kappadelta \|n_k\|$, so, from above,
\bequation
  \baligned
    v(x_k) - v(x_k + s_k) \geq \thalf\kappavm \kappadelta^{-2} \lambdav_k \|s_k\|^2 - \thalf \Hv_{Lip} \|s_k\|^3.
  \ealigned
\eequation
  From here, by Steps~\ref{step.phase1.f-iter-updates} and \ref{step.phase1.v-iter-rho} of Algorithm~\ref{alg.phase1} and under \eqref{eq.largeenoughlambda-1}, the result follows.
\eproof

The next lemma reveals upper and lower bounds for an important ratio that will hold during the iteration immediately following a \Viteration{} contraction.

\blemma\label{lem.lmdv/n_bound-1}
  For all $k \in \Ccal^v$, it follows that
  \bequation \label{eq.lmdv/n_bound-1}
    \sigmamin \leq \frac{\lambdav_{k+1}}{\|n_{k+1}\|} \leq \max \left \{\sigmamax, \(\frac{\gammalambda}{\gammac}\) \frac{\lambdav_k}{\|n_k\|} \right \}.
  \eequation
\elemma
\bproof
  The result follows using the same logic as the proof of \cite[Lemma 3.17]{CurtRobiSama16}.  
%
\iftechreport
  In particular, there are four cases to consider.

  \textbf{Case 1.} Suppose that Step~\ref{step.regularized1-Viter} of Algorithm~\ref{alg.viter} is reached.  Then, $\deltav_{k+1} = \|n_{k+1}\| = \|n(\lambdav)\|$ and $\lambdav_{k+1} = \lambdav$, where $(n(\lambdav),\lambdav)$ is computed in Steps~\ref{step.lambdahat-Viter}--\ref{step.find-n} of Algorithm~\ref{alg.viter}. As Step~\ref{step.regularized1-Viter} of Algorithm~\ref{alg.viter} is reached, the condition in Step~\ref{step.checkratio-Viter} of Algorithm~\ref{alg.viter} holds. Therefore, $\lambdav_{k+1} / \|n_{k+1}\| \leq \sigmamax$. To find a lower-bound on the ratio, let $\Hv_k = V_k \Xi^v_k V_k^T$ where $V_k$ is an orthonormal matrix of eigenvectors and $\Xiv_k = \diag(\xiv_{k,1},\dots,\xiv_{k,n})$ with $\xiv_{k,1} \leq \dots \leq \xiv_{k,n}$ is a diagonal matrix of eigenvalues of $\Hv_k$. Since $k \in \Ccal^v \subseteq \Ical$, $\|\gv_k\| \geq \epsilon > 0$; therefore, $\lambdav = \hat{\lambda}^v > \lambdav_k$, leading to $\Hv_k + \lambdav I \succ 0$.  Thus,
  \bequationn
    \|n(\lambdav)\|^2 = \|V_k (\Xiv_k + \lambdav I)^{-1} V_k^T \gv_k\|^2 = {\gv_k}^T V_k (\Xiv_k + \lambdav I)^{-2} V_k^T g_k.
  \eequationn
  From orthonormality of $V_k$ and Steps~\ref{step.lambdahat-Viter}--\ref{step.find-n} of Algorithm~\ref{alg.viter}, it follows that
  \bequationn
    \frac{\|n(\lambdav)\|^2}{\|\gv_k\|^2} = \frac{{\gv_k}^T V_k (\Xiv_k + \lambdav I)^{-2} V_k^T g_k}{\|V_k^T \gv_k\|^2} \leq \(\xiv_{k,1}+\lambdav_k + (\sigmamin\|\gv_k\|)^{1/2}\)^{-2}.
  \eequationn
  Hence, since $\lambdav_k \geq \max\{0,-\xiv_{k,1}\}$, one finds
  \bequationn
    \frac{\lambdav_{k+1}}{\|n_{k+1}\|} = \frac{\lambdav}{\|n(\lambdav)\|} \geq \frac{(\lambdav_k + (\sigmamin\|\gv_k\|)^{1/2})(\xiv_{k,1} + \lambdav_k + (\sigmamin\|\gv_k\|)^{1/2})}{\|\gv_k\|} \geq \sigmamin.
  \eequationn
  
  \textbf{Case 2.}  Suppose that Step~\ref{step.deltamagic-Viter} of Algorithm~\ref{alg.viter} is reached.  Then, $\deltav_{k+1} = \|n_{k+1}\| = \|n(\lambdav)\|$ where $(n(\lambdav),\lambdav)$ is computed in Step~\ref{step.linsysmagic-Viter} of Algorithm~\ref{alg.viter}, meaning that
  \bequationn
    \frac{\lambdav_{k+1}}{\|n_{k+1}\|} = \frac{\lambdav}{\|n(\lambdav)\|}\ \ \text{where}\ \ \sigmamin \leq \frac{\lambdav}{\|n(\lambdav)\|} \leq \sigmamax.
  \eequationn
  
  The other two cases that may occur correspond to situations in which the condition in Step~\ref{step.lamcheck-Viter} of Algorithm~\ref{alg.viter} tests false, in which case $\lambdav_k > 0$ and the pair $(n(\lambdav),\lambdav)$ is computed in Steps~\ref{step.lambda-Viter}--\ref{step.linsys-Viter} of Algorithm~\ref{alg.viter}. This means, in particular, that
  \bequation \label{eq.ratiolowerbound1-1}
    \sigmamin \leq \frac{\lambdav_k}{\|n_k\|} \leq \frac{\lambdav}{\|n(\lambdav)\|},
  \eequation
  where the latter inequality follows since $\lambdav = \gammalambda \lambdav_k > \lambdav_k$ , which, in turn, implies by standard trust region theory that $\|n(\lambdav)\| < \|n_k\|$.  Let us now consider the cases.
  
  \textbf{Case 3.} Suppose that Step~\ref{step.delta2-Viter} of Algorithm~\ref{alg.viter} is reached.  Then, $\lambdav_{k+1} = \lambdav$ and $\|n_{k+1}\| = \deltav_{k+1}$.  In conjunction with~\eqref{eq.ratiolowerbound1-1}, it follows that
  \bequationn
    \sigmamin \leq \frac{\lambdav_{k+1}}{\|n_{k+1}\|} = \frac{\lambdav}{\|n(\lambdav)\|} = \frac{\gammalambda \lambdav_k}{\|n(\lambdav)\|} \leq \frac{\gammalambda \lambdav_k}{\gammac \|n_k\|},
  \eequationn
  where the last inequality follows from the condition in Step~\ref{eq.sabotage} in Algorithm~\ref{alg.viter}.
  
  \textbf{Case 4.} Suppose that Step~\ref{step.deltanormal-Viter} of Algorithm~\ref{alg.viter} is reached, so $\deltav_{k+1} = \gammac \|n_k\|$.  According to standard trust region theory, since $\|n(\lambdav)\| < \gammac \|n_k\|$, one can conclude that $\lambdav_k \leq \lambdav_{k+1} \leq \lambdav = \gammalambda \lambdav_k$. Hence, with~\eqref{eq.ratiolowerbound1-1}, it follows that
  \bequationn
    \sigmamin < \frac{\sigmamin}{\gammac} \leq \frac{\lambdav_k}{\gammac \|n_k\|} \leq \frac{\lambdav_{k+1}}{\|n_{k+1}\|} \leq \frac{\gammalambda \lambdav_k}{\gammac \|n_k\|}.
  \eequationn

  The result follows since we have obtained the desired inequalities in all cases.
\fi
\eproof

Now, we prove that the sequence $\{\sigmav_k\}$ is bounded.

\blemma\label{lem.sigmavbound-1}
  There exists $\sigmav_{\max} \in (0,\infty)$ such that $\sigmav_k \leq \sigmav_{\max}$ for all $k\in\Ical$.
\elemma
\bproof
  If $k \notin \Ccal^v$, then Steps~\ref{step.phase1.f-iter_sigmav} and \ref{step.phase1.v-iter_sigmav} of Algorithm~\ref{alg.phase1} give $\sigmav_{k+1} \gets \sigmav_k$.  Otherwise, if $k \in \Ccal^v$, meaning that $\rho_k^v < \kapparho$, then there are two cases to consider.  If $k \in \Ccal^v$ and $\lambdav_k < \sigmamin \|n_k\|$, then, by Steps~\ref{step.checkratio-Viter}--\ref{step.deltamagic-Viter} of Algorithm~\ref{alg.viter}, Step~\ref{step.phase1.v-iter_sigmav} of Algorithm~\ref{alg.phase1}, and the fact that $\lambdav_{k+1} = \lambdav$ and $n_{k+1} = n(\lambdav)$, where $(n(\lambdav,\lambdav)$ are computed either in Steps~\ref{step.lambda0-Viter}--\ref{step.find-n} or Step~\ref{step.linsysmagic-Viter} of Algorithm~\ref{alg.viter}, it follows that $\sigmav_{k+1} \leq \max\{\sigmav_k,\sigmamax\}$.  Finally, if $k \in \Ccal^v$ and $\lambdav_k \geq \sigmamin \|n_k\|$, then it follows from Lemma \ref{lem.lambdak-lowerbound-1} that
  \bequationn
    \lambdav_k < \kappadelta^2\kappavm^{-1}(\Hv_{Lip} + 2\kapparho)\|s_k\|.
  \eequationn
  From the fact that $\lambdav_k \geq \sigmamin \|n_k\|$, Lemma~\ref{lem.nonzero_steps}(i), \eqref{eq.kkt-normal.comp}, and \eqref{def.delta_s}, it follows that $\|s_k\| \leq \deltas_k \leq \kappadelta \deltav_k = \kappadelta \|n_k\|$. Hence, by~Step \ref{step.phase1.v-iter_sigmav} of Algorithm~\ref{alg.phase1} and Lemma~\ref{lem.lmdv/n_bound-1}, one finds
  \bequationn
    \sigmav_{k+1} \gets \max \left\{\sigmav_k, \frac{\lambdav_{k+1}}{\|n_{k+1}\|} \right\} 
    \leq \max \left\{\sigmav_k, \sigmamax, \(\frac{\gammalambda}{\gammac}\) \frac{\kappadelta^3}{\kappavm}(\Hv_{Lip} + 2\kapparho)\right\}.
  \eequationn
  Combining the results of these cases gives the desired conclusion.
\eproof

We now give a lower-bound for the norm of some types of successful steps.

\blemma \label{lem.gv-and-s-1}
  For all $k \in \Scal^v_{\sigma} \cup \Scal^f$, the accepted step $s_k$ satisfies
  \bequation
    \|s_k\| \geq (\Hv_{Lip} + \kappaht + \sigmav_{\max}/\kappantn^2)^{-1/2} \|\gv_{k+1}\|^{1/2}. 
  \eequation 
\elemma
\bproof
  Let $k \in \Scal^v_{\sigma} \cup \Scal^f$. It follows from \eqref{eq.kkt-normal.dual}, the Mean Value Theorem, the fact that $s_k = n_k + t_k$, Assumption~\ref{ass.M-lip}, Lemma~\ref{lem.nonzero_steps}(ii), and \eqref{eq.Ht_vs_s} that there exists a vector $\overline{x}$ on the line segment $[x_k,x_k+s_k]$ such that
  \begin{align}
    \|\gv_{k+1}\| &= \|\gv_{k+1} - \gv_k - (\Hv_k + I \lambdav_k) n_k \| \nonumber \\
    &= \|(\Hv(\overline{x}) - \Hv_k) s_k +  \Hv_k t_k - \lambdav_k n_k\| \nonumber \\
   &\leq \Hv_{Lip}\|s_k\|^2 +  \kappaht \|s_k\|^2 +  \frac{\lambdav_k}{\|n_k\|} \|n_k\|^2. \label{eq.g-and-s-1}
  \end{align}
  From Step~\ref{step.accept-Viter} of Algorithm~\ref{alg.viter} (if $k\in\Scal^v_\sigma$) and~\eqref{eq.f-iter_sigma_v} (if $k\in\Scal^f$), one finds $\lambdav_k/\|n_k\| \leq \sigmav_k$.  Combining this with~\eqref{eq.g-and-s-1}, Lemma~\ref{lem.nonzero_steps}(ii), \eqref{eq.s_vs_n}, and Lemma~\ref{lem.sigmavbound-1}, it follows that
  \bequationn
    \|\gv_{k+1}\| \leq \Hv_{Lip}\|s_k\|^2 + \kappaht \|s_k\|^2 + \sigmav_k \|n_k\|^2 \leq (\Hv_{Lip} + \kappaht + \sigmav_{\max}/\kappantn^2) \|s_k\|^2,
  \eequationn
  which gives the desired result.
\eproof

We now give an iteration complexity result for a subset of successful iterations.

\blemma\label{lem.numAccepted-Vsteps-2}
  For any $\epsilon \in (0,\infty)$, the total number of elements in
  \bequationn
    \Kcal(\epsilon) := \{k \in \Ical : k \geq 0\ \text{and}\ (k-1) \in \Scal^v_{\sigma} \cup \Scal^f\}
  \eequationn
  is at most
  \bequation \label{eq.numAccepted-Vsteps-2}
    \left \lfloor \(\frac{\vmax_0}{\kapparho(1-\kappavtwo)(\Hv_{Lip} + \kappaht + \sigmav_{\max}/\kappantn^2)^{-3/2}}\) \epsilon ^{-3/2} \right \rfloor =: K_\sigma(\epsilon) \geq 0.
  \eequation
\elemma
\bproof
  From Lemma~\ref{lem.vmax_decrease} and Lemma~\ref{lem.gv-and-s-1}, it follows that, for all $k \in \Kcal(\epsilon) \subseteq \Ical$,
\bequationn
  \baligned
    \vmax_{k-1} - \vmax_k &\geq \kapparho(1-\kappavtwo) \|s_{k-1}\|^3 \\
    &\geq \kapparho(1-\kappavtwo) (\Hv_{Lip} + \kappaht + \sigmav_{\max}/\kappantn^2)^{-3/2} \|\gv_k\|^{3/2} \\
    &\geq \kapparho(1-\kappavtwo) (\Hv_{Lip} + \kappaht + \sigmav_{\max}/\kappantn^2)^{-3/2} \epsilon^{3/2}.
  \ealigned
\eequationn
  In addition, since $|\Kcal(\epsilon)| < \infty$ follows by Theorem~\ref{thm.finitetermination}, the reduction in $\vmax_k$ obtained up to the largest index in $\Kcal(\epsilon)$, call it $\kbar(\epsilon)$, satisfies
  \bequalin
    \vmax_0 - \vmax_{\kbar(\epsilon)} &= \sum_{k = 1}^{\kbar_{\epsilon}} (\vmax_{k-1} - \vmax_k) \geq \sum_{k \in \Kcal(\epsilon)} (\vmax_{k-1} - \vmax_k) \\
    &\geq |\Kcal(\epsilon)| \kapparho(1-\kappavtwo) (\Hv_{Lip} + \kappaht + \sigmav_{\max}/\kappantn^2)^{-3/2} \epsilon^{3/2}.
  \eequalin
  Rearranging this inequality to yield an upper bound for $|\Kcal(\epsilon)|$ and using the fact that $\vmax_k \geq 0$ for all $k \in \Ical$ (see Lemma~\ref{lem.vmax}), the desired result follows.
\eproof

In order to bound the total number of successful iterations in $\Ical$, we also need an upper bound for the cardinality of $\Scal^v_{\deltamax}$.  This is the subject of our next lemma.

\blemma\label{lem.numSvdeltamax}
  The cardinality of the set $\Scal^v_{\deltamax}$ is bounded above by
  \bequation \label{eq.numSvdeltamax}
    \left \lfloor \frac{\vmax_0}{\kapparho \kappantn^3(1-\kappavtwo){(\deltamaxv_0)}^3} \right \rfloor := K_{\deltamax}^v \geq 0.
  \eequation
\elemma
\bproof
  For all $k \in \Scal^v_{\deltamax} \subseteq \Scal$, it follows from Lemma~\ref{lem.vmax_decrease},  Lemma~\ref{lem.nonzero_steps}(ii), \eqref{eq.s_vs_n}, and Lemma~\ref{lem.deltav-changes}(ii) that the decrease in the trust funnel radius satisfies
  \bequationn
    \baligned
      \vmax_k - \vmax_{k+1} &\geq \kapparho(1-\kappavtwo)\|s_k\|^3 \geq \kapparho \kappantn^3(1-\kappavtwo)\|n_k\|^3 \\
      &= \kapparho \kappantn^3(1-\kappavtwo){(\deltamaxv_k)}^3 \geq \kapparho \kappantn^3(1-\kappavtwo){(\deltamaxv_0)}^3.
    \ealigned
  \eequationn
  Now, using the fact that $\{\vmax_k\}$ is bounded below by zero (see Lemma~\ref{lem.vmax}), one finds
  \bequationn
    \vmax_0 \geq \sum_{k \in \Scal^v_{\deltamax}} (\vmax_k - \vmax_{k+1}) \geq |\Scal^v_{\deltamax}| \kapparho \kappantn^3(1-\kappavtwo){(\deltamaxv_0)}^3,
  \eequationn
  which gives the desired result.
\eproof

Having now provided upper bounds for the numbers of successful iterations, we need to bound the number of unsuccessful iterations in $\Ical$.  To this end, first we prove that a critical ratio increases by at least a constant factor after an iteration in $\Ccal^v$.

\blemma\label{lem.ratio-increase-viter1-1}
  If $k \in \Ccal^v$ and $\lambdav_k \geq \sigmamin \|n_k\|$, then
  \bequationn
    \frac{\lambdav_{k+1}}{\|n_{k+1}\|} \geq \min \left \{\gammalambda, \frac{1}{\gammac} \right \} \(\frac{\lambdav_k}{\|n_k\|}\).
  \eequationn
\elemma
\bproof
  The proof follows the same logic as in \cite[Lemma 3.23]{CurtRobiSama16}. 
%
\iftechreport
In particular, since $k \in \Ccal^v$ and, with Lemma~\ref{lem.nonzero_steps}(i), it follows that $\lambdav_k \geq \sigmamin \|n_k\| > 0$, one finds that the condition in Step~\ref{step.lamcheck-Viter} of Algorithm~\ref{alg.viter} tests false. Hence, $(n(\lambdav),\lambdav)$ is computed in Steps~\ref{step.lambda-Viter}--\ref{step.linsys-Viter} of Algorithm~\ref{alg.viter} so that $\lambdav = \gammalambda \lambdav_k > \lambdav_k$ and $n(\lambdav)$ solves $\mathcal{Q}_k^v(\lambdav)$.  Let us now consider the two cases that may occur.
  
  \textbf{Case 1.} Suppose that Step~\ref{step.delta2-Viter} of Algorithm~\ref{alg.viter} is reached, meaning that $\|n(\lambdav)\| \geq \gammac \|n_k\|$. It follows that $\|n_{k+1}\| = \deltav_{k+1} < \deltav_k = \|n_k\|$ and $\lambdav_{k+1} = \gammalambda \lambdav_k$, i.e.,
  \bequation
    \frac{\lambdav_{k+1}}{\|n_{k+1}\|} > \frac{\gammalambda \lambdav_k}{\|n_k\|}.
  \eequation
  
  \textbf{Case 2.} Suppose that Step~\ref{step.deltanormal-Viter} of Algorithm~\ref{alg.viter} is reached, meaning that $\|n(\lambdav)\| < \gammac \|n_k\|$. It follows that $\|n_{k+1}\| = \deltav_{k+1} = \gammac \|n_k\|$ and $\lambdav_{k+1} \geq \lambdav_k$. Consequently,
  \bequation
    \frac{\lambdav_{k+1}}{\|n_{k+1}\|} \geq \frac{\lambdav_k}{\gammac \|n_k\|}.
  \eequation
  The result now follows from the conclusions of these two cases.
\fi
\eproof

We are now able to provide an upper bound on the number of unsuccessful iterations in $\Ccal^v$ that may occur between any two successful iterations.

\blemma\label{lem.numofVcontract}
  If $\kbar \in \Scal \cup \{0\}$, then
  \bequation
    |\Ccal^v \cap \Ical_{\Scal}(\kbar)| \leq 1 + \left\lfloor \frac{1}{\log(\min\{\gammalambda,\gammac^{-1}\})}\log\(\frac{\sigmav_{\max}}{\sigmamin}\)\right\rfloor =: K_{\Ccal}^v \geq 0.
  \eequation
\elemma
\bproof
  The result holds trivially if $|\Ccal^v \cap \Ical_{\Scal}(\kbar)| = 0$. Thus, we may proceed under the assumption that $|\Ccal^v \cap \Ical_{\Scal}(\kbar)| \geq 1$.  Let $k_{\Ccal^v}$ be the smallest element in $\Ccal^v \cap \Ical_{\Scal}(\kbar)$.  It then follows from Lemma~\ref{lem.deltav-changes}(i)-(ii), Lemma~\ref{lem.not-expansion}, and Step~\ref{step.phase1.f-iter-updates} of Algorithm~\ref{alg.phase1} that for all $k \in \Ical$ satisfying $k_{\Ccal^v} + 1 \leq k \leq k_{\Scal}(\kbar)$ we have
  \bequationn
    \|n_k\| \leq \deltav_k \leq  \deltav_{k_{\Ccal^v}+1} < \deltav_{k_{\Ccal^v}} \leq \deltamaxv_{k_{\Ccal^v}} \leq \deltamaxv_{k_{\Scal}(\kbar)},
  \eequationn
  which for $k = k_{\Scal}(\kbar)$ means that $k_{\Scal}(\kbar) \in \Scal^f \cup \Scal^v_{\sigma}$.  From Lemma~\ref{lem.lmdv/n_bound-1}, it follows that $\lambdav_{k_{\Ccal^v}+1} \geq \sigmamin \|n_{k_{\Ccal^v}+1}\|$, which by $k_{\Scal}(\kbar) \in \Scal^f \cup \Scal^v_{\sigma}$, Lemma~\ref{lem.sigmavbound-1}, Lemma~\ref{lem.ratio-increase-viter1-1}, \eqref{eq.f-iter_sigma_v}, Step~\ref{step.accept-Viter} of Algorithm~\ref{alg.viter}, and the fact that $(n_{k+1},\lambdav_{k+1}) = (n_{k},\lambdav_{k})$ for any $k \in \Ccal^f$ means
  \bequationn
    \sigmav_{\max} \geq \sigmav_{k_{\Scal}(\kbar)} \geq \frac{\lambdav_{k_{\Scal}(\kbar)}}{\|n_{k_{\Scal}(\kbar)}\|} \geq \(\min\left\{\gammalambda,\frac{1}{\gammac}\right\}\)^{|\Ccal^v \cap \Ical_{\Scal}(\kbar)|-1} \sigmamin,
  \eequationn
  from which the desired result follows.
\eproof

For our ultimate complexity result, the main component that remains to prove is a bound on the number of unsuccessful iterations in $\Ccal^f$ between any two successful iterations.  To this end, we first need some preliminary results pertaining to the trial step and related quantities during an \Fiteration.  Our first such result pertains to the change in the objective function model yielded by the tangential step.

\blemma \label{lem.n->s-fdecrease}
  For any $k \in \Ical$, the vectors $n_k$ and $t_k$ and dual variable $\lambdaf_k$ satisfy
  \bequation \label{eq.n->s-fdecrease}
    \mf_k(n_k) - \mf_k(n_k + t_k) = \thalf t_k^T (H_k + \lambdaf_k I) t_k + \thalf \lambdaf_k \|t_k\|^2 + \lambdaf_k n_k^T t_k.
  \eequation
\elemma
\bproof
  If $k\notin\Ical^t$ so that $t_k = 0$ and $\lambdaf_k = 0$ (by the \computesteps\ subroutine in Algorithm~\ref{alg.phase1}), then~\eqref{eq.n->s-fdecrease} trivially holds.  Thus, for the remainder of the proof, let us assume that $k\in\Ical^t$.  It now follows from the definition of $\mf_k$ that
  \bequalin
     &\ \mf_k(n_k) - \mf_k(n_k + t_k) \\
    =&\ g_k^T n_k + \thalf n_k^T H_k n_k - g_k^T (n_k+t_k) - \thalf (n_k+t_k)^T H_k (n_k+t_k) \\
    =&\ - (g_k + H_k n_k)^T t_k - \thalf t_k^T H_k t_k \\
    =&\ - (g_k + (H_k +\lambdaf_k I)n_k + (H_k + \lambdaf_k I) t_k + J_k^T \yf_k)^T t_k \\
     &\ + \thalf t_k^T (H_k + \lambdaf_k I) t_k + \thalf \lambdaf_k \|t_k\|^2 + \lambdaf_k n_k^T t_k + {(\yf_k)}^T J_k t_k \\
    =&\ \thalf t_k^T (H_k + \lambdaf_k I) t_k + \thalf \lambdaf_k \|t_k\|^2 + \lambdaf_k n_k^T t_k,
  \eequalin
  where the last equality follows from \eqref{eq.kkt_tangent.dual}.
\eproof

The next lemma reveals that, for an \Fiteration{}, if the dual variable for the tangential step trust region constraint is large enough, then the trust region constraint is active and the iteration will be successful.

\blemma\label{lem.lambdaf-lowerbound}
  For all $k \in \Fcal$, if the trial step $s_k$ and the dual variable $\lambdaf_k$ satisfy
  \bequation \label{eq.lambdat-lowerbound}
    \lambdaf_k \geq (\kappafm\kappast^2 (1-\kappantt))^{-1} (\kappahs + H_{Lip} + 2\kapparho) \|s_k\|,
  \eequation
  then $\|s_k\| = \deltas_k$ and $\rho_k^f \geq \kapparho$.  
\elemma
\bproof
  Observe from \eqref{eq.lambdat-lowerbound} and Lemma~\ref{lem.nonzero_steps}(i) that $\lambdaf_k > 0$, which along with~\eqref{eq.kkt_tangent.comp} proves that $\|s_k\| = \deltas_k$.  Next, since $k \in \Fcal$, it must mean that \eqref{eq.f-iteration-conditions} is satisfied.  It then follows from \eqref{eq.f-iter_mf_red_t}, Lemma~\ref{lem.n->s-fdecrease}, \eqref{eq.kkt_tangent}, \eqref{eq.f-iter_nt_vs_t}, and~\eqref{eq.f-iter_t_vs_s} that
  \begin{align}
    \mf_k(0) - \mf_k(s_k)
      &\geq \kappafm (\mf_k(n_k) - \mf_k(s_k)) \nonumber \\
      &= \kappafm(\thalf t_k^T (H_k + \lambdaf_k I) t_k + \thalf \lambdaf_k \|t_k\|^2 + \lambdaf_k n_k^T t_k) \nonumber \\
      &\geq \kappafm (\thalf-\thalf\kappantt) \lambdaf_k \|t_k\|^2 \geq \thalf \kappafm \kappast^2 (1 - \kappantt) \lambdaf_k \|s_k\|^2. \label{eq.diff.m}
  \end{align}
  Next, the Mean Value Theorem gives the existence of an $\xbar \in [x_k, x_k+s_k]$ such that $f(x_k + s_k) = f_k + g_k^T s_k + \half s_k^T \nabla^2 f(\xbar) s_k$, which with \eqref{eq.f-iter_Hs_vs_s} and Assumption~\ref{ass.M-lip} gives
  \bequalin
        &\ \mf_k(s_k) - f(x_k + s_k) \\
    =   &\ f_k + g_k^T s_k + \thalf s_k^T H_k s_k - f_k - g_k^T s_k - \thalf s_k^T \nabla^2 f(\xbar) s_k \\
    =   &\ \thalf \(s_k^TH_ks_k - s_k^T\nabla^2 f(x_k)s_k + s_k^T\nabla^2 f(x_k)s_k - s_k^T\nabla^2 f(\xbar)s_k\) \\
    =   &\ \thalf s_k^T\(H_k - \nabla^2 f(x_k)\)s_k + \thalf s_k^T\(\nabla^2 f(x_k) - \nabla^2 f(\xbar)\)s_k \\
    \geq&\ -\thalf \left\|\(H_k - \nabla^2 f(x_k)\)s_k \right\| \|s_k\| - \thalf \left\|\(\nabla^2 f(x_k) - \nabla^2 f(\xbar)\)s_k \right\| \|s_k\| \\  
    \geq&\ -\thalf \(\kappahs + H_{Lip}\) \|s_k\|^3. 
  \eequalin
  Finally, combining the previous inequality, $f_k = \mf_k(0)$, and~\eqref{eq.diff.m}, one finds
  \bequalin
    f_k - f(x_k + s_k)
      &=    f_k - \mf_k(s_k) + \mf_k(s_k) - f(x_k + s_k) \\
      &\geq \thalf \kappafm \kappast^2 (1-\kappantt) \lambdaf_k \|s_k\|^2 - \thalf (\kappahs + H_{Lip}) \|s_k\|^3,
  \eequalin
  which combined with \eqref{eq.lambdat-lowerbound} shows that $\rho_k^f \geq \kapparho$ as desired.
\eproof

We now show that a critical ratio increases by at least a constant factor after any unsuccessful \Fiteration\ followed by an iteration in which a nonzero tangential step is computed and not reset to zero.

\blemma\label{lem.ratio-increase-Fiter1}
  If $k \in \Ccal^f$, $\lambdaf_k \geq \sigmamin \|s_k\|$, and $(k+1) \in \Ical^t$, then
  \bequationn
    \frac{\lambdaf_{k+1}}{\|s_{k+1}\|} \geq \(\frac{1}{\gammac}\)\frac{\lambdaf_k}{\|s_k\|}.
  \eequationn
\elemma
\bproof
  With Lemma~\ref{lem.nonzero_steps}(i), it follows that $\lambdaf_k \geq \sigmamin \|s_k\| > 0$, meaning that $\|s_k\| = \deltas_k$.  In addition, since $k \in \Ccal^f$, one finds that the condition in Step~\ref{step.lamcheck-Fiter} of Algorithm~\ref{alg.Fiter} tests false in iteration $k$.  Hence, Step~\ref{step.deltanormal-Fiter} of Algorithm~\ref{alg.Fiter} is reached, meaning, with \eqref{def.delta_s}, that $\deltaf_{k+1} = \gammac \|s_k\| \leq \gammac \kappadelta \deltav_k$.  Then, from the facts that $\gammac < 1$ and $\deltav_{k+1} \gets \deltav_k$ (see Step~\ref{step.phase1.f-iter-updates} of Algorithm~\ref{alg.phase1}), it follows that $\deltaf_{k+1} \leq \kappadelta \deltav_{k+1}$.  Consequently, again with~\eqref{def.delta_s}, it follows that $\|s_{k+1}\| = \deltas_{k+1} = \deltaf_{k+1} = \gammac \|s_k\|$.  Combining this with the fact that Lemma~\ref{lem.deltaf-changes}(i) yields $\lambdaf_{k+1} \geq \lambdaf_k$, the result follows.
\eproof

\blemma\label{lem.lmdf/s_bound}
  If $k \in \Ccal^f$ and $(k+1) \in \Ical^t$, then $\sigmamin \leq \lambdaf_{k+1}/\|s_{k+1}\|$.
\elemma
\bproof
  Since $k \in \Ccal^f$, there are two cases to consider.
  
  \textbf{Case 1:} Step~\ref{step.deltamagic-Fiter} of Algorithm~\ref{alg.Fiter} is reached.  In this case, it follows that $\|s_{k+1}\| = \deltaf_{k+1} = \|n_k+t(\lambdaf)\|$ with $(t(\lambdaf),\lambdaf)$ computed in Step~\ref{step.linsysmagic-Fiter} of Algorithm~\ref{alg.Fiter}.  Together with the fact that $(k+1)\in\Ical^t$, it follows that $\lambdaf_{k+1}/\|s_{k+1}\| = \lambdaf/\|n_k+t(\lambdaf)\| \geq \sigmamin$.

  \textbf{Case 2:} Step~\ref{step.deltamagic-Fiter} of Algorithm~\ref{alg.Fiter} is not reached.  This only happens if the condition in Step~\ref{step.lamcheck-Fiter} of Algorithm~\ref{alg.Fiter} tested false, meaning that $\lambdaf_k / \|s_k\| \geq \sigmamin$. Hence, from Lemma~\ref{lem.ratio-increase-Fiter1}, it follows that $\lambdaf_{k+1}/\|s_{k+1}\| \geq \lambdaf_k/\gammac \|s_k\|$, which by the facts that $\gammac < 1$ and $\lambdaf_k / \|s_k\| \geq \sigmamin$ gives the desired result.
\eproof

Next, we provide a bound on the number of iterations in $\Ccal^f$ that may occur before the first or between consecutive iterations in the set $\Scal\cup\Vcal$.

\blemma\label{lem.numofFcontract}
  If $\kbar \in \Scal \cup \Vcal \cup \{0\}$, then
  \bequation
    |\Ical_{\Scal \cup \Vcal}(\kbar)| \leq 2 + \left\lfloor \frac{1}{\log(\gammac^{-1})}\log\(\frac{\kappahs + H_{Lip} + 2\kapparho}{\sigmamin \kappafm \kappast^2 (1-\kappantt)}\)\right\rfloor =: K_{\Ccal}^f \geq 0.
  \eequation
\elemma
\bproof
  Let $\kbar \in \Scal \cup \Vcal \cup \{0\}$.  Then, $\Ical_{\Scal\cup\Vcal}(\kbar) \subseteq \Ccal^f$.  The result follows trivially if $|\Ical_{\Scal \cup \Vcal}(\kbar)| \leq 1$.  Therefore, for the remainder of the proof, let us assume that $|\Ical_{\Scal \cup \Vcal}(\kbar)| \geq 2$.  It follows from Lemma~\ref{lem.lmdf/s_bound}, $\kbar+1 \in \Ccal^f$, and $\kbar+2 \in \Ccal^f \subseteq \Fcal$ (meaning that $t_{\kbar+2} \neq 0$ and $(\kbar+2)\in\Ical^t$) that $\sigmamin \leq \lambdaf_{\kbar+2}/\|s_{\kbar+2}\|$.  Combining this inequality with Lemma~\ref{lem.lambdaf-lowerbound}, 
Lemma~\ref{lem.ratio-increase-Fiter1}, the fact that and $(k_{\Scal\cup\Vcal}(\kbar)-1) \in \Ccal^f$ to get
  \bequalin
    \sigmamin\(\frac{1}{\gammac}\)^{(k_{\Scal\cup\Vcal}(\kbar)-1)-(\kbar+2)}
     \leq \frac{\lambdaf_{k_{\Scal\cup\Vcal}(\kbar)-1}}{\|s_{k_{\Scal\cup\Vcal}(\kbar)-1}\|} 
    \leq \(\frac{\kappahs + H_{Lip} + 2\kapparho}{\kappafm \kappast^2 (1-\kappantt)} \).
  \eequalin
  The desired result now follows since $|\Ical_{\Scal \cup \Vcal}(\kbar)| = k_{\Scal\cup\Vcal}(\kbar)-\kbar-1$.
\eproof

We have now arrived at our complexity result for phase 1.

\btheorem \label{thm.complexity}
  For a scalar $\epsilon \in (0,\infty)$, the cardinality of $\Ical$ is at most
  \bequation\label{eq.K_epsilon}
    K(\epsilon) := 1 + (K_\sigma(\epsilon) + K_{\deltamax}^v)(K_{\Ccal}^v + 1)K_{\Ccal}^f, 
  \eequation
  where $K_\sigma(\epsilon)$, $K_\Delta^v$, $K_\Ccal^v$, and $K_\Ccal^f$ are defined in Lemmas~\ref{lem.numAccepted-Vsteps-2}, \ref{lem.numSvdeltamax}, \ref{lem.numofVcontract}, and \ref{lem.numofFcontract}, respectively.  Consequently, for any $\bar\epsilon \in (0,\infty)$, it follows that $K(\epsilon) = \Ocal(\epsilon^{-3/2})$ for all $\epsilon \in (0,\bar\epsilon)$.
\etheorem
\bproof
  Without loss of generality, let us assume that at least one iteration is performed.  Then, Lemmas~\ref{lem.numAccepted-Vsteps-2} and \ref{lem.numSvdeltamax} guarantee that at most $K_\sigma(\epsilon) + K_\Delta^v$ successful iterations are included in $\Ical$.  In addition, Lemmas~\ref{lem.1-expansion}, \ref{lem.numofVcontract}, and \ref{lem.numofFcontract} guarantee that, before each successful iteration, there can be at most $(K_{\Ccal}^v + 1)K_{\Ccal}^f$ unsuccessful iterations.  Also accounting for the first iteration, the desired result follows.
\eproof

If the constraint Jacobians encountered by the algorithm are not rank deficient (and do not tend toward rank deficiency), then the following corollary gives a similar result as that above, but for an infeasibility measure.

\bcorollary\label{cor.complexity}
  Suppose that, for all $k \in \N{}$, the constraint Jacobian $J_k$ has full row rank with singular values bounded below by $\zeta_{\min} \in (0,\infty)$.  Then, for $\epsilon \in (0,\infty)$, the cardinality of $\Ical_c := \{k \in \N{} : \|c_k\| > \epsilon/\zeta_{\min}\}$, is at most $K(\epsilon)$ defined in \eqref{eq.K_epsilon}.  Consequently, for any $\bar\epsilon \in (0,\infty)$, the cardinality of $\Ical_c$ is $\Ocal(\epsilon^{-3/2})$ for all $\epsilon \in (0,\bar\epsilon)$.
\ecorollary
\bproof
  Under the stated conditions, $\|\gv_k\| \equiv \|J_k^T c_k\| \geq \zeta_{\min} \|c_k\|$ for all $k \in \Ical$.  Thus, since $\|\gv_k\| \leq \epsilon$ implies $\|c_k\| \leq \epsilon/\zeta_{\min}$, the result follows from~Theorem~\ref{thm.complexity}.
\eproof

\section{Phase 2: Obtaining Optimality}\label{sec.phase2}

A complete algorithm for solving problem~\eqref{prob.opt} proceeds as follows.  The phase 1 method, Algorithm~\ref{alg.phase1}, is run until either an approximate feasible point or approximate infeasible stationary point is found, i.e., for some $(\epsilon_{feas},\epsilon_{inf}) \in (0,\infty) \times (0,\infty)$, the method is run until, for some $k \in \N{}$,
\bsubequations\label{eq.near_feas_or_infeas}
  \begin{align}
    \|c_k\| &\leq \epsilon_{feas} \label{eq.near_feas} \\ \text{or}\ \ \|J_k^Tc_k\| &\leq \epsilon_{inf}\|c_k\|. \label{eq.near_infeas}
  \end{align}
\esubequations
If phase 1 terminates with \eqref{eq.near_feas} failing to hold and \eqref{eq.near_infeas} holding, then the entire algorithm is terminated with a declaration of having found an infeasible (approximately) stationary point.  Otherwise, if \eqref{eq.near_feas} holds, then a phase 2 method is run that maintains at least $\epsilon_{feas}$-feasibility while seeking optimality.

With this idea in mind, how should the termination tolerance $\epsilon$ in Algorithm~\ref{alg.phase1} be set so that \eqref{eq.near_feas_or_infeas} is achieved within at most $\Ocal(\epsilon^{-3/2})$ iterations, as is guaranteed by the analysis in the previous section?  Given $(\epsilon_{feas},\epsilon_{inf}) \in (0,\infty) \times (0,\infty)$, we claim that Algorithm~\ref{alg.phase1} should be employed with $\epsilon = \epsilon_{feas}\epsilon_{inf}$.  Indeed, with this choice, if the final point produced by phase 1, call it~$x_k$, does not yield~\eqref{eq.near_feas}, then it must satisfy $\|\gv_k\| \equiv \|J_k^Tc_k\|/\|c_k\| \leq \epsilon/\epsilon_{feas} = \epsilon_{inf}$, which is exactly~\eqref{eq.near_infeas}.


There are various options for phase 2, three of which are worth mentioning.
\bitemize
  \item Respecting the current state-of-the-art nonlinear optimization methods, one can run a trust funnel method such as that in \cite{GoulToin10}.  One can even run such a method with the initial trust funnel radius for $v(x) = \thalf\|c(x)\|^2$ set at $\thalf\epsilon_{feas}^2$ so that $\epsilon_{feas}$-feasibility will be maintained as optimality is sought.  We do not claim worst-case iteration complexity guarantees for such a method, though empirical evidence suggests that such a method would perform well.  This is the type of approach for which experimental results are provided in \S\ref{sec.numerical}.
  \item With an eye toward attaining good complexity properties, one can run the objective-target-following approach proposed as \cite[Alg.~4.1,~Phase~2]{CartGoulToin13a}.  This approach essentially applies an \ARC{} algorithm for unconstrained optimization \cite{CartGoulToin11a,CartGoulToin11b} (see also the previous work in \cite{Grie81,NestPoly06,WeisDeufErdm07}) to minimize the residual function $\Phi : \R{N} \times \R{} \to \R{}$ defined by $\Phi(x,t) = \|c(x)\|^2 + \|f(x) - t\|^2$.  In iteration $k \in \N{}$, the subsequent iterate $x_{k+1}$ is computed to reduce $\Phi(\cdot,t_k)$ as in \ARC{} while the subsequent target $t_{k+1}$ is chosen to ensure, amongst other relationships, that $t_{k+1} \leq t_k$ and $|f_k - t_k| \leq \epsilon_{feas}$ for all $k \in \N{}$, where it is assumed that $\epsilon_{feas} \in (0,1)$.  In \cite{CartGoulToin13a}, it is shown that, for the phase 2 algorithm with $\epsilon \in (0,\epsilon_{feas}^{1/3}]$, the number of iterations required to generate a primal iterate $x_k$ satisfying \eqref{eq.near_feas} and either the \emph{relative} KKT error condition
  \bequationn
    \|g_k + J_k^Ty_k\| \leq \epsilon\|(y_k,1)\|\ \ \text{for some}\ \ y_k \in \R{M}
  \eequationn
  or the constraint violation stationarity condition
  \bequationn
    \|J_k^Tc_k\| \leq \epsilon\|c_k\|
  \eequationn
  is at most $\Ocal(\epsilon^{-3/2}\epsilon_{feas}^{-1/2})$.  This should be viewed in two ways.  First, if $\epsilon = \epsilon_{feas}^{2/3}$, then the overall complexity is $\Ocal(\epsilon_{feas}^{-3/2})$, though of course this corresponds to a looser tolerance on the relative KKT error than on feasibility.  Second, if $\epsilon = \epsilon_{feas}$ (so that the two tolerances are equal), then the overall complexity is $\Ocal(\epsilon_{feas}^{-2})$.  We claim that an approach based on \TRACE{} \cite{CurtRobiSama16} (instead of \ARC{}) could instead be employed yielding the same worst-case iteration complexity properties; see \iftechreport Appendix~\ref{app.phase2}\else \cite{CurtRobiSama17}\fi.
  \item Finally, let us point out that in cases that $c$ is affine, one could run an optimization method, such as the \ARC{} method from \cite{CartGoulToin11a,CartGoulToin11b} or the \TRACE{} method from \cite{CurtRobiSama16}, where steps toward reducing the objective function are restricted to the null space of the constraint Jacobian.  For such a reduced-space method, $\epsilon_{feas}$-feasibility will be maintained while the analyses in \cite{CartGoulToin11a,CartGoulToin11b,CurtRobiSama16} guarantee that the number of iterations required to reduce the norm of the reduced gradient below a given tolerance $\epsilon_{opt} \in (0,\infty)$ is at most $\Ocal(\epsilon_{opt}^{-3/2})$.  With $\epsilon = \epsilon_{opt} = \epsilon_{feas}$, this gives an overall (phase 1 $+$ phase 2) complexity of $\Ocal(\epsilon^{-3/2})$, which matches the optimal complexity for the unconstrained case.
\eitemize

\section{Numerical Experiments} \label{sec.numerical}

Our goal in this section is to demonstrate that instead of having a phase~1 method that solely seeks (approximate) feasibility (such as in \cite{CartGoulToin13a}), it is beneficial to employ a phase~1 method such as ours that simultaneously attempts to reduce the objective function.  To show this, a Matlab implementation of our phase~1 method, Algorithm~\ref{alg.phase1}, has been written.  The implementation has two modes: one following the procedures of Algorithm~\ref{alg.phase1} and one employing the same procedures except that the tangential step $t_k$ is set to zero for all $k \in \N{}$ so that all iterations are \Viteration{}s.  We refer to the former implementation as \TrustFunnel\ and the latter as \TrustVonly.  For phase~2 for both methods, following the current state-of-the-art, we implemented a trust funnel method based on that proposed in~\cite{GoulToin10} with the modification that the normal step computation is never skipped.  In both phases~1 and~2, all subproblems are solved to high accuracy using a Matlab implementation of the trust region subproblem solver described as~\cite[Alg.~7.3.4]{ConGT00a}, which in large part goes back to the work in \cite{MoreSore83}.  The fact that the normal step computation is never skipped and the subproblems are always solved to high accuracy allows our implementation to ignore so-called ``y-iterations''~\cite{GoulToin10}.

Phase~1 in each implementation terminates in iteration $k \in \N{}$ if either
\bequation\label{eq.chunky}
  \|c_k\|_\infty \leq 10^{-6}\max\{\|c_0\|_\infty,1\}\ \ \text{or}\ \ \left\{ \baligned \|J_k^Tc_k\|_\infty &\leq 10^{-6}\max\{\|J_0^Tc_0\|_\infty,1\} \\ \text{and}\ \|c_k\|_\infty &> 10^{-3}\max\{\|c_0\|_\infty,1\} \ealigned \right.
\eequation
whereas Phase~2 terminates in iteration $k \in \N{}$ if either the latter pair of conditions above holds or, with $y_k$ computed as least squares multipliers for all $k \in \N{}$, if
\bequationn
  \|g_k + J_k^Ty_k\|_\infty \leq 10^{-6}\max\{\|g_0 + J_0^Ty_0\|_\infty,1\}.
\eequationn
Input parameters used in the code are stated in Table~\ref{tab.parameters}.  The only values that do not appear are $\kapparho$ and $\gammac$.  For $\kapparho$, for simplicity we employed this constant in \eqref{eq.f-iter_vmax_rule} and \eqref{eq.vmax_update_f} as well as in the step acceptance conditions in Step~\ref{step.accept-Fiter} in Algorithm~\ref{alg.Fiter} and Step~\ref{step.accept-Viter} in Algorithm~\ref{alg.viter}.  That said, our convergence analysis is easily adapted to handle different values in these places: our code uses $\kapparho = 10^{-12}$ in \eqref{eq.f-iter_vmax_rule} and  \eqref{eq.vmax_update_f} but $\kapparho = 10^{-8}$ in the step acceptance conditions.  For $\gammac$, our code uses $0.5$ in the context of an \Fiteration{} (Algorithm~\ref{alg.Fiter}) and $10^{-2}$ in the context of a \Viteration{} (Algorithm~\ref{alg.viter}), where again our analysis easily allows using different constants in these places.

\btable[ht]
  \centering
  \caption{Input parameters for \TrustFunnel\ and \TrustVonly.}
  \label{tab.parameters}
  \texttt{
  \btabular{|l|l||l|l||l|l||l|l|}
    \hline
    $\kappan$   & 9e-01 & $\kappast$   & 1e-12     & $\kappap$   & 1e-06 & $\kappadelta$  & 1e+02 \\
    $\kappavm$  & 1e-12 & $\kappantt$  & 1-(2e-12) & $\kappaht$  & 1e+20 & $\gammae$      & 2e+00 \\
    $\kappantn$ & 1e-12 & $\kappavone$ & 9e-01     & $\kappahs$  & 1e+20 & $\gammalambda$ & 2e+00 \\
    $\kappafm$  & 1e-12 & $\kappavtwo$ & 9e-01     & $\sigmamin$ & 1e-12 & $\sigmamax$    & 1e+20 \\
    \hline
  \etabular
  }
\etable

We ran \TrustFunnel\ and \TrustVonly\ to solve the equality constrained problems in the CUTEst test set \cite{GoulOrbaToin13}.  Among 190 such problems, we removed 78 that had a constant (or null) objective, 13 for which phase 1 of both algorithms terminated immediately at the initial point due to the former condition in \eqref{eq.chunky}, three for which both algorithms terminated phase 1 due to the latter pair of conditions in \eqref{eq.chunky} (in each case within one iteration), two on which both algorithms encountered a function evaluation error,
and one on which both algorithms failed due to small stepsizes (less than $10^{-20}$) in phase 1.  We also removed all problems on which neither algorithm terminated within one hour.  The remaining set consisted of 33 problems.

\btable
  \centering
  \tiny
  \caption{Numerical results for \TrustFunnel\ and \TrustVonly.}
  \label{tab.numericalresults}
  \rowcolors{1}{white}{lightgray}
  \renewcommand{\tabcolsep}{0.18cm}
  \texttt{
  \btabular{|l|r|r|r|r|r|r|r|r|r|r|r|r|r|}
    \hline
    &  &  & \multicolumn{6}{c|}{\TrustFunnel} & \multicolumn{5}{c|}{\TrustVonly} \\ \cline{4-14}
    &  &  & \multicolumn{4}{c|}{Phase 1} & \multicolumn{2}{c|}{Phase 2} & \multicolumn{3}{c|}{Phase 1} & \multicolumn{2}{c|}{Phase 2} \\ \cline{4-14}
    \multicolumn{1}{|c|}{Problem} & \multicolumn{1}{c|}{$n$} & \multicolumn{1}{c|}{$m$} & \multicolumn{1}{c|}{\#V} & \multicolumn{1}{c|}{\#F} & \multicolumn{1}{c|}{$f$} & \multicolumn{1}{c|}{$\|g+J^Ty\|$} & \multicolumn{1}{c|}{\#V} & \multicolumn{1}{c|}{\#F} & \multicolumn{1}{c|}{\#V} & \multicolumn{1}{c|}{$f$} & \multicolumn{1}{c|}{$\|g+J^Ty\|$} & \multicolumn{1}{c|}{\#V} & \multicolumn{1}{c|}{\#F} \\ \hline
           BT1 &       2 &      1 &      4 &      0 & -8.02e-01 & +4.79e-01 &      0 &    139 &       4 & -8.00e-01 & +7.04e-01 &      7 &    136  \\ \hline 
      BT10 &       2 &      2 &     10 &      0 & -1.00e+00 & +5.39e-04 &      1 &      0 &      10 & -1.00e+00 & +6.74e-05 &      1 &      0  \\ \hline 
      BT11 &       5 &      3 &      6 &      1 & +8.25e-01 & +4.84e-03 &      2 &      0 &       1 & +4.55e+04 & +2.57e+04 &     16 &     36  \\ \hline 
      BT12 &       5 &      3 &     12 &      1 & +6.19e+00 & +1.18e-05 &      0 &      0 &      16 & +3.34e+01 & +4.15e+00 &      4 &      8  \\ \hline 
       BT2 &       3 &      1 &     22 &      8 & +1.45e+03 & +3.30e+02 &      3 &     12 &      21 & +6.14e+04 & +1.82e+04 &      0 &     40  \\ \hline 
       BT3 &       5 &      3 &      1 &      0 & +4.09e+00 & +6.43e+02 &      1 &      0 &       1 & +1.01e+05 & +8.89e+02 &      0 &      1  \\ \hline 
       BT4 &       3 &      2 &      1 &      0 & -1.86e+01 & +1.00e+01 &     20 &     12 &       1 & -1.86e+01 & +1.00e+01 &     20 &     12  \\ \hline 
       BT5 &       3 &      2 &     15 &      2 & +9.62e+02 & +2.80e+00 &     14 &      2 &       8 & +9.62e+02 & +3.83e-01 &      3 &      1  \\ \hline 
       BT6 &       5 &      2 &     11 &     45 & +2.77e-01 & +4.64e-02 &      1 &      0 &      14 & +5.81e+02 & +4.50e+02 &      5 &     59  \\ \hline 
       BT7 &       5 &      3 &     15 &      6 & +1.31e+01 & +5.57e+00 &      5 &      1 &      12 & +1.81e+01 & +1.02e+01 &     19 &     28  \\ \hline 
       BT8 &       5 &      2 &     50 &     26 & +1.00e+00 & +7.64e-04 &      1 &      1 &      10 & +2.00e+00 & +2.00e+00 &      1 &     97  \\ \hline 
       BT9 &       4 &      2 &     11 &      1 & -1.00e+00 & +8.56e-05 &      1 &      0 &      10 & -9.69e-01 & +2.26e-01 &      5 &      1  \\ \hline 
  BYRDSPHR &       3 &      2 &     29 &      2 & -4.68e+00 & +1.28e-05 &      0 &      0 &      19 & -5.00e-01 & +1.00e+00 &     16 &      5  \\ \hline 
     CHAIN &     800 &    401 &      9 &      0 & +5.12e+00 & +2.35e-04 &      3 &     20 &       9 & +5.12e+00 & +2.35e-04 &      3 &     20  \\ \hline 
       FLT &       2 &      2 &     15 &      4 & +2.68e+10 & +3.28e+05 &      0 &     13 &      19 & +2.68e+10 & +3.28e+05 &      0 &     17  \\ \hline 
   GENHS28 &      10 &      8 &      1 &      0 & +9.27e-01 & +5.88e+01 &      0 &      0 &       1 & +2.46e+03 & +9.95e+01 &      0 &      1  \\ \hline 
  HS100LNP &       7 &      2 &     16 &      2 & +6.89e+02 & +1.74e+01 &      4 &      1 &       5 & +7.08e+02 & +1.93e+01 &     14 &      3  \\ \hline 
  HS111LNP &      10 &      3 &      9 &      1 & -4.78e+01 & +4.91e-06 &      2 &      0 &      10 & -4.62e+01 & +7.49e-01 &     10 &      1  \\ \hline 
      HS27 &       3 &      1 &      2 &      0 & +8.77e+01 & +2.03e+02 &      3 &      5 &       1 & +2.54e+01 & +1.41e+02 &     11 &     34  \\ \hline 
      HS39 &       4 &      2 &     11 &      1 & -1.00e+00 & +8.56e-05 &      1 &      0 &      10 & -9.69e-01 & +2.26e-01 &      5 &      1  \\ \hline 
      HS40 &       4 &      3 &      4 &      0 & -2.50e-01 & +1.95e-06 &      0 &      0 &       3 & -2.49e-01 & +3.35e-02 &      2 &      1  \\ \hline 
      HS42 &       4 &      2 &      4 &      1 & +1.39e+01 & +3.94e-04 &      1 &      0 &       1 & +1.50e+01 & +2.00e+00 &      3 &      1  \\ \hline 
      HS52 &       5 &      3 &      1 &      0 & +5.33e+00 & +1.54e+02 &      1 &      0 &       1 & +8.07e+03 & +4.09e+02 &      0 &      1  \\ \hline 
       HS6 &       2 &      1 &      1 &      0 & +4.84e+00 & +1.56e+00 &     32 &    136 &       1 & +4.84e+00 & +1.56e+00 &     32 &    136  \\ \hline 
       HS7 &       2 &      1 &      7 &      1 & -2.35e-01 & +1.18e+00 &      7 &      2 &       8 & +3.79e-01 & +1.07e+00 &      5 &      2  \\ \hline 
      HS77 &       5 &      2 &     13 &     30 & +2.42e-01 & +1.26e-02 &      0 &      0 &      17 & +5.52e+02 & +4.54e+02 &      3 &     11  \\ \hline 
      HS78 &       5 &      3 &      6 &      0 & -2.92e+00 & +3.65e-04 &      1 &      0 &      10 & -1.79e+00 & +1.77e+00 &      2 &     30  \\ \hline 
      HS79 &       5 &      3 &     13 &     21 & +7.88e-02 & +5.51e-02 &      0 &      2 &      10 & +9.70e+01 & +1.21e+02 &      0 &     24  \\ \hline 
   MARATOS &       2 &      1 &      4 &      0 & -1.00e+00 & +8.59e-05 &      1 &      0 &       3 & -9.96e-01 & +9.02e-02 &      2 &      1  \\ \hline 
      MSS3 &    2070 &   1981 &     12 &      0 & -4.99e+01 & +2.51e-01 &     50 &      0 &      12 & -4.99e+01 & +2.51e-01 &     50 &      0  \\ \hline 
   MWRIGHT &       5 &      3 &     17 &      6 & +2.31e+01 & +5.78e-05 &      1 &      0 &       7 & +5.07e+01 & +1.04e+01 &     12 &     20  \\ \hline 
  ORTHREGB &      27 &      6 &     10 &     15 & +7.02e-05 & +4.23e-04 &      0 &      6 &      10 & +2.73e+00 & +1.60e+00 &      0 &     10  \\ \hline 
   SPIN2OP &     102 &    100 &     57 &     18 & +2.04e-08 & +2.74e-04 &      0 &      1 &    time & +1.67e+01 & +3.03e-01 &   time &   time  \\ \hline 
  \etabular
  }
\etable

The results we obtained are provided in Table~\ref{tab.numericalresults}.  For each problem, we indicate the number of variables ($n$), number of equality constraints ($m$), number of \Viteration{}s~(\texttt{\#V}), number of \Fiteration{}s~(\texttt{\#F}), objective function value at the end of phase~1~($f$), and dual infeasibility value at the end of phase~1 ($\|g+J^Ty\|$).  We use \texttt{time} for \texttt{SPIN2OP} for \TrustVonly\ to indicate that it hit the one hour time limit (after 350 phase~1 iterations) without terminating.  The results illustrate that, within a comparable number of iterations, our trust funnel algorithm, represented by~\TrustFunnel, typically yields \emph{better} final points from phase~1.  This can be seen in the fact that the objective at the end of phase~1, dual infeasibility at the end of phase~1, and the number of iterations required in phase~2 are all typically smaller for \TrustFunnel\ than they are for \TrustVonly.  Note that for some problems, such as \texttt{BT1}, \TrustFunnel\ only performs \Viteration{}s in phase~1, yet yields a better final point than does \TrustVonly; this occurs since the phase~1 iterations in \TrustFunnel\ may involve nonzero tangential steps.

\section{Conclusion}\label{sec.conclusion}

An algorithm has been proposed for solving equality constrained optimization problems.  Following the work in \cite{CartGoulToin13a}, but based on trust funnel and trust region ideas from \cite{CurtRobiSama16,GoulToin10}, the algorithm represents a next step toward the design of practical methods for solving constrained optimization problems that offer strong worst-case iteration complexity properties.  In particular, the algorithm involves two phases, the first seeking (approximate) feasibility and the second seeking optimality, where a key contribution is the fact that improvement in the objective function is sought in both phases.  If a phase~2 method such as that proposed in \cite{CartGoulToin13a} is employed, then the overall algorithm attains the same complexity properties as the method in~\cite{CartGoulToin13a}.  The results of numerical experiments show that the proposed method benefits by respecting the objective function in both phases.

\bibliographystyle{plain}
\bibliography{references}

\iftechreport

\appendix

\section{Phase 2 Details}\label{app.phase2}

The goal of this appendix is to show that a phase~2 method can be built upon the \TRACE{} algorithm from \cite{CurtRobiSama16} yielding the same worst-case iteration complexity properties as the \ARC-based method in \cite{CartGoulToin13a}.  We state a complete phase~2 algorithm, then prove similar properties for it as those proved for the method in \cite{CartGoulToin13a}.  For the algorithm and our analysis of it, we make the following additional assumption.

\bassumption\label{ass.f-bounded}
  For all $x \in \R{N}$ with $\|c(x)\| \leq \epsilon_{feas} \in (0,\infty)$, the objective function $f$ is bounded from below and above by $f_{\min} \in \R{}$ and $f_{\max} \in \R{}$, respectively.  In addition, the problem functions $f$ and $c$ and their first and second derivatives are Lipschitz continuous on the path defined by all phase~2 iterates.
\eassumption

As previously mentioned, the phase~2 method is in many ways based on applying an algorithm for solving unconstrained optimization problems to minimize the residual function $\Phi(x,t) = \thalf \|r(x,t)\|^2$ where $r : \R{N} \times \R{} \to \R{}$ is defined by
\bequation \label{eq.def-residual}
  r(x,t) = \bpmatrix c(x) \\ f(x) - t \epmatrix.
\eequation
Updated dynamically by the algorithm, the parameter $t$ may be viewed as a target value for reducing the objective function value.

The phase 2 algorithm is stated as Algorithm~\ref{alg.nftrace}.  We refer the reader to \cite{CartGoulToin13a} for further details on the design of the algorithm and to \cite{CurtRobiSama16} for further details on \TRACE{}.

\balgorithm[ht]
  \small
  \caption{\TRACE{} Algorithm for Phase 2}
  \label{alg.nftrace}
  \balgorithmic[1]
    \smallskip
    \Require termination tolerance $\epsilon \in (0,\infty)$ and $x_0 \in \R{N}$ with $\|c_0\| \leq \epsilon_{feas} \in (0,\infty)$
    \smallskip
    \AlgBreak
    \Procedure{\TRACE\_\textsc{phase}\_\textsc{2}}{}
      \State set $t_0 \gets f_0 - \sqrt{\epsilon_{feas}^2 - \|c_0\|^2}$ \label{step.t0}
      \For{$k \in \N{}$}
        \State perform one iteration of \TRACE{} toward minimizing $\Phi(x,t_k)$ to compute $s_k$
        \If{$s_k$ is an acceptable step}
          \State set $x_{k+1} \gets x_k + s_k$ (and other quantities following \TRACE{})
          \If{$r(x_{k+1},t_k) \neq 0$ and $\|\nabla_x \Phi(x_{k+1},t_k)\| \leq \epsilon \|r(x_{k+1},t_k)\|$} \label{step.checktermination}
            \State terminate
        \Else
          \State set $t_{k+1} \gets f(x_{k+1}) - \sqrt{\|r(x_k,t_k)\|^2 - \|r(x_{k+1},t_k)\|^2 + \(f(x_{k+1}) - t_k\)^2}$ \label{step.t-update}
        \EndIf
      \Else
        \State set $x_{k+1} \gets x_k$ (and other quantities following \TRACE{})
        \State set $t_{k+1} \gets t_k$
      \EndIf
      \EndFor
    \EndProcedure
    \smallskip
  \ealgorithmic
\ealgorithm

The following lemma, whose proof follows that of~\cite[Lemma 4.1]{CartGoulToin13a}, states some useful properties of the generated sequences.

\blemma \label{lem.decreasingt}
  For all $k \in \N{}$, it follows that
  \bsubequations\label{eq.phase2_basic}
    \begin{align}
      t_{k+1} &\leq t_k, \label{eq.t-nonincreasing} \\
      0 \leq f(x_k) - t_k &\leq \epsilon_{feas}, \label{eq.fbigerthant} \\
      \|r(x_k,t_k)\| &= \epsilon_{feas}, \label{eq.r-epsilon} \\
      \text{and}\ \ \|c(x_k)\| &\leq \epsilon_{feas}. \label{eq.stayfeasible}
    \end{align}
  \esubequations
\elemma
\bproof
  Note that, in \TRACE{}, the objective function is monotonically nonincreasing; see \cite[Eq.~(2.5);~Alg.~1,~Step~5]{CurtRobiSama16}.  Hence, each acceptable step $s_k$ computed in Algorithm~\ref{alg.nftrace} yields $\Phi(x_{k+1},t_k) \leq \Phi(x_k,t_k)$, from which it follows that the value for $t_{k+1}$ in Step~\ref{step.t-update} is well-defined.  Then, since all definitions and procedures in Algorithm~\ref{alg.nftrace} that yield \eqref{eq.phase2_basic} are exactly the same as in \cite[Alg.~4.1]{CartGoulToin13a}, a proof for the inequalities in \eqref{eq.phase2_basic} is given by the proof of~\cite[Lemma 4.1]{CartGoulToin13a}.
\eproof

In the next lemma, we recall a critical result from \cite{CurtRobiSama16}, arguing that it remains true for Algorithm~\ref{alg.nftrace} under our assumptions about the problem functions.

\blemma \label{lem.sigmabound_phase2}
  Let $\{\sigma_k\}$ be generated as in \TRACE{} \cite{CurtRobiSama16}.  Then, there exists a scalar constant $\sigma_{\max} \in (0,\infty)$ such that $\sigma_k \leq \sigma_{\max}$ for all $k \in \N{}$. 
\elemma
\bproof
  The result follows in a similar manner as \cite[Lem.~3.18]{CurtRobiSama16}.  Here, similar to \cite[\S5]{CartGoulToin13a}, it is important to note that Assumption~\ref{ass.f-bounded} ensures that $\Phi$ and its first and second derivatives are globally Lipschitz continuous on a path defined by the phase~2 iterates.  This ensures that results of the kind given as \cite[Lem.~3.16--3.17]{CurtRobiSama16} hold true, which are necessary for proving \cite[Lem.~3.18]{CurtRobiSama16}.
\eproof

We now argue that the number of iterations taken for any fixed value of the target for the objective function is bounded above by a positive constant.

\blemma \label{lem.fixed-t-iterations}
  The number of iterations required before the first accepted step or between two successive accepted steps with a fixed target $t$ is bounded above by
  \bequationn
    K_t := 2 +  \left \lfloor \frac{1}{\log(\min\{\gammalambda,\gammac^{-1}\})} \log \(\frac{\sigma_{\max}}{\sigmamin}\) \right \rfloor,
  \eequationn
where the constants $\gammalambda\in(0,\infty)$, $\gammac\in(0,1)$, are $\sigmamin\in(0,\infty)$ are parameters used by \TRACE{} (see~\cite[Alg.~1]{CurtRobiSama16}) that are independent of $k$ and satisfy $\sigmamin \leq \sigma_{\max}$.  
\elemma
\bproof
  The properties of \TRACE{} corresponding to so-called \emph{contraction} and \emph{expansion} iterations all hold for Algorithm~\ref{alg.nftrace} for sequences of iterations in which a target value is held fixed.  Therefore, the result follows by \cite[Lem.~3.22 and Lem.~3.24]{CurtRobiSama16}, which combined show that the maximum number of iterations of interest is equal to the maximum number of contractions that may occur plus one.
\eproof

The next lemma merely states a fundamental property of \TRACE{}.

\blemma \label{lem.functiondecrease}
  Let $H_\Phi \in (0,\infty)$ be the Lipschitz constant for the Hessian function of $\Phi$ 
  along the path of phase~2 iterates and let $\eta \in (0,1)$ be the acceptance constant from \TRACE{}.  Then, for $x_{k+1}$ following an accepted step $s_k$, it follows that
  \bequationn
    \Phi(x_k, t_k) - \Phi(x_{k+1},t_k) \geq \eta (H_\Phi + \sigma_{\max})^{-3/2} \|\nabla_x \Phi(x_{k+1},t_k)\|^{3/2}.
  \eequationn
\elemma
\bproof
  With Lemma~\ref{lem.sigmabound_phase2} and adapting the conclusion of~\cite[Lem.~3.19]{CurtRobiSama16}, the result follows as in the beginning of the proof of \cite[Lem.~3.20]{CurtRobiSama16}.
\eproof

The preceding lemma allows us to prove the following useful result.

\blemma \label{lem.r-change}
  For $x_{k+1}$ following an accepted step $s_k$ yielding
  \bequation\label{eq.phiphiphi}
    \|\nabla_x \Phi(x_{k+1},t_k)\| > \epsilon\|r(x_{k+1},t_k)\|
  \eequation
  with $\epsilon$ the constant used in Algorithm~\ref{alg.nftrace}, it follows that
  \bequationn
    \|r(x_k,t_k)\| - \|r(x_{k+1},t_k)\| \geq \kappa_t \min\{\epsilon^{3/2} \epsilon_{feas}^{1/2}, \epsilon_{feas}\},
  \eequationn
  where $\beta \in (0,1)$ is any fixed problem-independent constant, $\omega := \eta (H_\Phi + \sigma_{\max})^{-3/2} \in (0,\infty)$ is the constant appearing in Lemma~\ref{lem.functiondecrease}, and $\kappa_t := \min\{\omega \beta^{3/2} , 1-\beta \}$.
\elemma
\bproof
  Along with the result of Lemma \ref{lem.functiondecrease}, it follows that
  \bequalin
    \|r(x_k,t_k)\|^2 - \|r(x_{k+1},t_k)\|^2 &\geq 2\omega \|\nabla_x \Phi(x_{k+1},t_k)\|^{3/2} \\
      &= 2\omega \(\frac{\|\nabla_x \Phi(x_{k+1},t_k)\|}{\|r(x_{k+1},t_k)\|}\)^{3/2} \|r(x_{k+1},t_k)\|^{3/2} \\
      &\geq 2\omega \epsilon^{3/2} \|r(x_{k+1},t_k)\|^{3/2}.
  \eequalin
  Then, if $\|r(x_{k+1},t_k)\| > \beta \|r(x_k,t_k)\|$, it follows with \eqref{eq.r-epsilon} that
  \bequation \label{eq.r2decrease}
    \|r(x_k,t_k)\|^2 - \|r(x_{k+1},t_k)\|^2 \geq 2\omega \epsilon^{3/2} \beta^{3/2} \|r(x_k,t_k)\|^{3/2} = 2\omega \epsilon^{3/2} \beta^{3/2} \epsilon_{feas}^{3/2},
  \eequation
  from which it follows along with $\|r(x_{k+1},t_k)\| \leq \|r(x_k,t_k)\|$ that
  \bequalin
    \|r(x_k,t_k)\| - \|r(x_{k+1},t_k)\| &= \frac{\|r(x_k,t_k)\|^2 - \|r(x_{k+1},t_k)\|^2}{\|r(x_k,t_k)\| + \|r(x_{k+1},t_k)\|} \\
      &\geq \frac{\|r(x_k,t_k)\|^2 - \|r(x_{k+1},t_k)\|^2}{2\|r(x_k,t_k)\|} \\
      &= \frac{\|r(x_k,t_k)\|^2 - \|r(x_{k+1},t_k)\|^2}{2\epsilon_{feas}} \geq \omega \epsilon^{3/2} \beta^{3/2} \epsilon_{feas}^{1/2}.
  \eequalin
  On the other hand, if $\|r(x_{k+1},t_k)\| \leq \beta \|r(x_k,t_k)\|$, then using~\eqref{eq.r-epsilon} it follows that
  \bequationn
    \|r(x_k,t_k)\| - \|r(x_{k+1},t_k)\| \geq (1-\beta) \|r(x_k,t_k)\| = (1-\beta) \epsilon_{feas}.
  \eequationn
  Combining the results of both cases, the desired conclusion follows.
\eproof

The following is an intermediate result used to prove the subsequent lemma.  We merely state the result for ease of reference;  see \cite[Lemma 5.2]{CartGoulToin13a} and its proof.

\blemma \label{lem.fromthere}
  Consider the following optimization problem in two variables:
  \bequationn
    \min_{(f,c) \in \R{2}}\ -f + \sqrt{\epsilon^2 - c^2}\ \ \st\ \ f^2 + c^2 \leq \tau^2,
  \eequationn
  where $0 < \tau < \epsilon$.  The global minimum of this problem is attained at $(f_*,c_*) = (\tau,0)$ with the optimal value given by $-\tau + \epsilon$.
\elemma

We next prove a lower bound on the decrease of the target value.

\blemma \label{lem.t-decrease-bound}
  Suppose that the termination tolerance is set so that $\epsilon \leq \epsilon_{feas}^{1/3}$.  Then, for $x_{k+1}$ following an accepted step $s_k$ such that the termination conditions in Step~\ref{step.checktermination} are not satisfied, it follows that, with $\kappa_t \in (0,1)$ defined as in Lemma~\ref{lem.r-change},
  \bequation\label{eq.ttt}
    t_k - t_{k+1} \geq \kappa_t \epsilon^{3/2} \epsilon_{feas}^{1/2}.
  \eequation
\elemma
\bproof
  A proof follows similarly to that of \cite[Lem.~5.3]{CartGoulToin13a}.  In particular, if the reason that the termination conditions in Step~\ref{step.checktermination} are not satisfied is because \eqref{eq.phiphiphi} holds, then Lemma~\ref{lem.r-change} and the fact that $\epsilon \leq \epsilon_{feas}^{1/3}$ imply that
  \bequationn
    \|r(x_k,t_k)\| - \|r(x_{k+1},t_k)\| \geq \kappa_t \min\{\epsilon^{3/2} \epsilon_{feas}^{1/2}, \epsilon_{feas}\} \geq \kappa_t \epsilon^{3/2} \epsilon_{feas}^{1/2}.
  \eequationn
  On the other hand, if the reason the termination conditions in Step~\ref{step.checktermination} are not satisfied is because $\|r(x_{k+1},t_k)\| = 0$, it follows from \eqref{eq.r-epsilon}, $\kappa_t \in (0,1)$, and $\epsilon \leq \epsilon_{feas}^{1/3}$ that
  \bequationn
    \|r(x_k,t_k)\| - \|r(x_{k+1},t_k)\| = \epsilon_{feas} \geq \kappa_t \epsilon^{3/2} \epsilon_{feas}^{1/2}.
  \eequationn
  Combining these two cases, \eqref{eq.def-residual}, and \eqref{eq.r-epsilon}, one finds that
  \bequationn
    (f(x_{k+1}) - t_k)^2 + \|c(x_{k+1})\|^2 = \|r(x_{k+1},t_k)\|^2 \leq (\epsilon_{feas} - \kappa_t\epsilon^{3/2}\epsilon_{feas}^{1/2})^2.
  \eequationn
  Now, from Step~\ref{step.t-update} of Algorithm~\ref{alg.nftrace}, \eqref{eq.def-residual}, and \eqref{eq.r-epsilon}, it follows that
  \bequalin
    t_k - t_{k+1} &= -(f(x_{k+1}) - t_k) + \sqrt{\|r(x_k,t_k)\|^2 - \|r(x_{k+1},t_k)\|^2 + \(f(x_{k+1}) - t_k\)^2} \\
      &= -(f(x_{k+1}) - t_k) + \sqrt{\|r(x_k,t_k)\|^2 - \|c(x_{k+1})\|^2} \\
      &= -(f(x_{k+1}) - t_k) + \sqrt{\epsilon_{feas}^2 - \|c(x_{k+1})\|^2}.
  \eequalin
  Overall, it follows that Lemma~\ref{lem.fromthere} can be applied (with ``$f$'' $= f(x_{k+1}) - t_k$, ``$c$'' $= \|c(x_{k+1})\|$, ``$\epsilon$'' $= \epsilon_{feas}$, and ``$\tau$'' $= \epsilon_{feas} - \kappa_t \epsilon^{3/2} \epsilon_{feas}^{1/2}$) to obtain the result.
\eproof

We now show that if the termination condition in Step~\ref{step.checktermination} of Algorithm~\ref{alg.nftrace} is never satisfied, then the algorithm takes infinitely many accepted steps.

\blemma \label{lem.A_infinite}
 If Algorithm~\ref{alg.nftrace} does not terminate finitely, then it takes infinitely many accepted steps.
\elemma
\bproof
  To derive a contradiction, suppose that the number of accepted steps is finite.  Then, since it does not terminate finitely, there exists $\kbar \in \N{}$ such that $s_k$ is not acceptable for all $k \geq \kbar$.  Therefore, by the construction of the algorithm, it follows that $t_k = t_{\kbar}$ for all $k \geq \kbar$.  This means that the algorithm proceeds as if the \TRACE{} algorithm from \cite{CurtRobiSama16} is being employed to minimize $\Phi(\cdot,t_{\kbar})$, from which it follows by \cite[Lemmas 3.7, 3.8, and 3.9]{CurtRobiSama16} that, for some sufficiently large $k \geq \kbar$, an acceptable step $s_k$ will be computed.  This contradiction completes the proof.
\eproof

Before proceeding, let us discuss the situation in which the termination conditions in Step~\ref{step.checktermination} of Algorithm~\ref{alg.nftrace} are satisfied. This discussion, originally presented in \cite{CartGoulToin13a}, justifies the use of these termination conditions.

Suppose $\|r(x_{k+1},t_k)\| \neq 0$ and $\|\nabla_x \Phi(x_{k+1},t_k)\| \leq \epsilon\|r(x_{k+1},t_k)\|$.  If $f_{k+1} = t_k$, then these mean that $\|c_{k+1}\| \neq 0$ and $\|\nabla_x \Phi(x_{k+1},t_k)\| \leq \epsilon\|c_{k+1}\|$, which along with $\nabla_x \Phi(x_{k+1},t_k) = J_{k+1}^T c_{k+1} + (f_{k+1} - t_k) g_{k+1}$ would imply that
\bequationn
  \frac{\|J_{k+1}^T c_{k+1}\|}{\|c_{k+1}\|} \leq \epsilon.
\eequationn
That is, under these conditions, $x_{k+1}$ is an approximate first-order stationary point for minimizing $\|c\|$. If $f_{k+1} \neq t_k$, then the satisfied termination conditions imply that
\bequationn
  \frac{\|J_{k+1}^T c_{k+1}+(f_{k+1} - t_k) g_{k+1}\|}{\|r(x_{k+1},t_k)\|} \leq \epsilon.
\eequationn
By dividing the numerator and denominator of the left-hand side of this inequality by $f(x_{k+1}) - t_k > 0$ (recall \eqref{eq.fbigerthant}), defining 
\begin{equation} \label{y-new-def}
y(x_{k+1},t_k) := c(x_{k+1})/(f(x_{k+1}) - t_k) \in \R{M},
\end{equation}
and substituting $y(x_{k+1},t_k)$ back into the inequality, one finds that
\bequation\label{eq.relativegradient2}
  \frac{\|J_{k+1}^T y(x_{k+1},t_k) + g_{k+1}\|}{\|(y(x_{k+1},t_k),1)\|} \leq \epsilon.
\eequation
As argued in \cite{CartGoulToin13a}, one may use a perturbation argument to say that $(x_{k+1},y(x_{k+1},t_k))$ satisfying the \emph{relative KKT error} conditions \eqref{eq.relativegradient2} and $\|c_{k+1}\| \leq \epsilon_{feas}$ corresponds to a first-order stationary point for problem~\eqref{prob.opt}.  Specifically, consider $x = x_*+\delta_x$ and $y = y_*+\delta_y$ where $(x_*,y_*)$ is a primal-dual pair satisfying the KKT conditions for problem~\eqref{prob.opt}. Then, a first-order Taylor expansion of $J(x_*)^T y_* + g(x_*)$ to estimate its value at $(x,y)$ yields the estimate $(\nabla^2 f(x_*) + \sum_{i=1}^M [y_i]_* \nabla^2 c_i(x^*)) \delta_x + J(x^*)^T \delta_y$. The presence of the dual variable $y^*$ in this estimate confirms that the magnitude of the dual variable should not be ignored in a relative KKT error condition such as \eqref{eq.relativegradient2}.

We now prove our worst-case iteration complexity result for phase 2.

\btheorem\label{lem.A_num}
Suppose that the termination tolerances are set so that $\epsilon \leq \epsilon_{feas}^{1/3}$ with $\epsilon_{feas} \in (0,1)$.  Then, Algorithm~\ref{alg.nftrace} requires at most $\Ocal{(\epsilon^{-3/2} \epsilon_{feas}^{-1/2})}$ iterations until 
the termination condition in Step~\ref{step.checktermination} is satisfied, at which point either
\bequation \label{eq.relativekkt}
  \frac{\|J_{k+1}^T y(x_{k+1},t_k) + g_{k+1}\|}{\|(y(x_{k+1},t_k),1)\|} \leq \epsilon
    \ \text{and} \ \|c_{k+1}\| \leq \epsilon_{feas}
  \eequation
  or 
  \bequation \label{eq.first-critical-c}
    \frac{\|J_{k+1}^T c_{k+1}\|}{\|c_{k+1}\|} \leq \epsilon
    \ \text{and} \ 
    \|c_{k+1}\| \leq \epsilon_{feas}
  \eequation
is satisfied, with $y(x_{k+1},t_k)$ defined in~\eqref{y-new-def}.   
\etheorem
\bproof
  Recall that if the termination condition in Step~\ref{step.checktermination} is satisfied for some $k \in \N{}$, then, by the arguments prior to the lemma,  either \eqref{eq.relativekkt} or \eqref{eq.first-critical-c} will be satisfied.  Thus, we aim to show an upper bound on the number of iterations required by the algorithm until the termination condition in Step~\ref{step.checktermination} is satisfied.
  
  Without loss of generality, let us suppose that the algorithm performs at least one iteration.  Then, we claim that there exists some $k \in \N{}$ such that the termination condition does not hold for $(x_k,t_{k-1})$, but does hold for $(x_{k+1},t_k)$.  To see this, suppose for contradiction that the termination condition is never satisfied.  Then, by Lemma~\ref{lem.t-decrease-bound}, it follows that for all $k \in \N{}$ such that $s_k$ is acceptable one finds that \eqref{eq.ttt} holds.  This, along with Lemma~\ref{lem.A_infinite}, implies that $\{t_k\} \searrow -\infty$.  However, this along with \eqref{eq.fbigerthant} implies that $\{f_k\} \searrow -\infty$, which contradicts Assumption~\ref{ass.f-bounded}.
  
  Now, since the termination condition is satisfied at $(x_{k+1},t_k)$, but not in the iteration before, it follows that $s_k$ must be an acceptable step.  Hence, from Assumption~\ref{ass.f-bounded}, \eqref{eq.fbigerthant}, Lemma~\ref{lem.t-decrease-bound}, Step~\ref{step.t0}, it follows that
  \bequalin
    f_{\min} \leq f_k \leq t_k + \epsilon_{feas}
    &\leq t_0 - K_{\Acal} \kappa_t \epsilon^{3/2} \epsilon_{feas}^{1/2} + \epsilon_{feas} \\
    &\leq f(x_0) - K_{\Acal} \kappa_t \epsilon^{3/2} \epsilon_{feas}^{1/2} + \epsilon_{feas},
  \eequalin
  where $K_{\Acal}$ is the number of accepted steps prior to iteration $(k+1)$.  Rearranging and since $\epsilon_{feas} \in (0,1)$, one finds along with Assumption~\ref{ass.f-bounded} that
  \bequation \label{eq.Acceptedcount}
    K_{\Acal} \leq \left \lceil \frac{f_{\max} - f_{\min} + 1}{\kappa_t \epsilon^{3/2} \epsilon_{feas}^{1/2}} \right \rceil.
  \eequation 
  From \eqref{eq.Acceptedcount} and Lemma~\ref{lem.fixed-t-iterations}, the desired result follows.
\eproof



\fi

\end{document}